\documentclass[12pt]{article}
\usepackage{amsmath,amsfonts,amssymb,theorem}

{\setlength{\labelwidth}{3.5cm}\setlength{\labelsep}{1.00em}
\setlength{\leftmargin}{5.00cm}\setlength{\rightmargin}{3.00cm}
{\setlength{\labelwidth}{3.5cm}\setlength{\labelsep}{1.5em}
\setlength{\leftmargin}{5.00cm}\setlength{\rightmargin}{3.00cm}
\addtolength{\hoffset}{0.25cm}
 \usepackage[matrix, arrow]{xy}

\DeclareMathOperator{\Supp}{Supp}

\DeclareMathOperator{\mld}{mld}

 \numberwithin{equation}{subsubsection}
 \numberwithin{footnote}{subsection}

 \newtheorem{cor}[subsubsection]{Corollary}
 \newtheorem{lem}[subsubsection]{Lemma}
 
 \newtheorem{thm}[subsubsection]{Theorem}
 \newtheorem{conj}[subsubsection]{Conjecture}

{
 \theorembodyfont{\rmfamily}
 \newtheorem{defn}[subsubsection]{Definition}
 \newtheorem{exa}[subsubsection]{Example}
 \newtheorem{exa-cr}[subsubsection]{Example--Construction}

 \newtheorem{rem}[subsubsection]{Remark}

 \newtheorem{recall}[subsubsection]{Recall}

}
 \newcommand{\qed}{\ifhmode\unskip\nobreak\fi\quad\ensuremath\square}
 \newenvironment{proof}{\paragraph{Proof}}{\par\medskip}

 \newcommand{\ke}[1]{$\acute{\mbox{e}}$}
 \newcommand{\ku}[1]{$\acute{\mbox{u}$}}
 \newcommand{\kl}[1]{$\acute{\mbox{l}}$}
 \newcommand{\kh}[1]{$\acute{\mbox{h}}$}
 \newcommand{\kr}[1]{$\acute{\mbox{r}}$}
 \newcommand{\kx}[1]{$\acute{\mbox{x}}$}
 \newcommand{\ki}[1]{${\^\i}$}

\title{ \large{\textbf{
Boundedness of $\epsilon$-lc Complements on Surfaces}}}

\author{\textbf{ Caucher Birkar{\thanks{Supported by the EPSRC and the Cecil King Prize.}}} \thanks{E-mail: caucher.birkar@maths.nott.ac.uk}}

\begin{document}

\maketitle

\tableofcontents

%%%%%%%%%%%%%%%%%%%%%%

\section{Boundedness of $\epsilon$-log canonical complements on surfaces}

\subsection{Introduction}

 The concept of {\textbf{complement}} was introduced and studied by Shokurov [Sh1, Sh2]. He used complements as a tool in the construction of 3-fold log flips [Sh1] and in the classification of singularities and contractions [Sh2]. Roughly speaking a complement is a ``good member'' of the anti-pluricanonical linear system i.e. a general member of $|-nK_X|$ for some $n>0$. In order, the existence of such member and the behaviour of the index $n$ are the most important problems in complement theory. Below we give the precise definition of the ``good member''. 

Throughout this paper we assume that the varieties involved are algebraic varieties over $\mathbb{C}$. In this section the varieties are all surfaces unless otherwise stated. By a log pair $(X/Z,B)$ we mean algebraic varieties $X$ and $Z$ equiped with a projective contraction $X\longrightarrow Z$ and $B$ is an $\mathbb{R}$-boundary on $X$. When we write $(X/P\in Z,B)$ we mean a log pair $(X/Z,B)$ with a fixed point $P\in Z$; in such situation we may shrink $Z$ around $P$ without mentioning it. The pair $(X/Z,B)$ is weak log Fano (WLF) if it has log canonical singularities (lc) and $-(K_X+B)$ is nef and big$/Z$ and $X$ is $\mathbb{Q}$-factorial. 

For the basic definitions of the Log Minimal Model Program (LMMP), the main references are [KMM] and [KM]. And to learn more about the complement theory [Sh2] and [Pr] are the best.

 \begin{defn}[Complements]\label{comp} Let $(X/Z,B=\sum_{i} b_{i}B_{i})$ be a $d$-dim pair where $X$ is normal and $B$ is an $\mathbb{R}$-boundary. Then $K_{X}+B^+$ is called an $(\epsilon, n)$-complement/$P\in Z$ (resp. in codim 2) for $K_X+B$, where $B^{+}=\sum_{i} b_{i}^{+}B_{i}$, if the following properties hold: 
 \begin{description}
 \item[$\diamond$] $(X,K_{X}+B^{+})$ is an $\epsilon$-lc pair/$P\in Z$ (resp. $\epsilon$-lc in codim 2) and $n(K_{X}+B^{+})\sim 0/P\in Z$.
 \item[$\diamond$] $\llcorner (n+1)b_{i}\lrcorner\leq nb_{i}^{+}$.
\end{description} 
 $K_{X}+B^+$ is called an $(\epsilon,\mathbb{R})$-complement/$P\in Z$ (resp.  in codim 2) for $K_X+B$ if $(X,K_{X}+B^{+})$   is $\epsilon$-lc$/P\in Z$ (resp. $\epsilon$-lc in codim 2), $K_{X}+B^{+}\sim_{\mathbb{R}} 0/P\in Z$ and $B^{+}\geq B$. An $(\epsilon,\mathbb{Q})$-complement/$P\in Z$ can be similarly defined where $\sim_{\mathbb{R}}$ is replaced by $\sim_{\mathbb{Q}}$.
\end{defn}

Despite the not quite easy definition above, complements have very good birational and inductive properties which make the theory a powerful tool to apply to the LMMP. In order, complements don't always exist even with strong conditions such as $-(K_X+B)$ being nef [Sh2, 1.1]. But for example they certainly do when $-(K_X+B)$ is nef and big and $K_X+B$ is lc. In this paper, in all the situations that occur the complements usually exist. So we concentrate on the second main problem about complements which is related to several open problems in the LMMP; namely the boundedness. We state conjectures on the boundedness of complements due to  Shokurov.

\begin{conj}[Weak $\epsilon$-lc Complements]\label{weak} 
Let $\Gamma \subseteq [0,1]$ be a set of real numbers which satisfies D.C.C. Then for any $0<\delta$ and $d$ there exist a finite set $\mathcal{N}_{\delta,d,\Gamma}$ of positive integers and $0<\epsilon$ such that any d-dim $\delta$-lc (resp. in codim 2) weak log Fano pair $(X/P\in Z,K_{X}+B)$, where $B\in \Gamma$, would be $(\epsilon,n)$-complementary/$P\in Z$  (resp. in codim 2) for some $n\in \mathcal{N}_{\delta,d ,\Gamma}$ .  
  
\end{conj}

We show the above conjecture as $WC_{\delta, d,\Gamma}$ for short. 

\begin{conj}[Strong $\epsilon$-lc Complements] \label{strong}
 For any $0<\epsilon$ and $d$ there exists a finite set $\mathcal{N}_{\epsilon, d}$ of positive integers such that any d-dim $\epsilon$-lc(resp. in codim 2) weak log Fano pair $(X/P\in Z,K_{X}+B)$ has an $(\epsilon,n)$-complement/$P\in Z$ (resp. in codim 2) for some $n\in \mathcal{N}_{\epsilon,d}$.    
\end{conj}

 We show the above conjecture as $SC_{\delta, d}$ for short. If we replace $\epsilon>0$ with $\epsilon=0$  in the conjecture above (it makes a big difference) then we get the usual conjecture on the boundedness of lc complements which has been studied by Shokurov, Prokhorov and others [Sh2, PSh, PSh1, Pr]. It is proved in dim 2 [Sh2] with some restrictions on the coefficients of $B$.

The following important conjecture, due to Alexeev and Borisov brothers, is related to the conjectures above [Mc, A1, PSh, MP]. 

\begin{conj}[BAB]\label{BAB}
Let $\delta >0$ be a real number and $\Gamma\subset [0,1]$. Then those varieties $X$ for which $(X,B)$ is a $\delta$-lc WLF pair of dim $d$ for a boundary $B\in \Gamma$, are elements of an algebraic family. 
\end{conj}

We show the above conjecture as $BAB_{\delta, d,\Gamma}$ for short. Alexeev proved  $BAB_{\delta, 2,\Gamma}$ for any $\delta>0$ and $\Gamma$ [A1]. This conjecture  was proved by Kawamata for terminal singularities in dim 3 [K1] and  $BAB_{1, 3,\{0\}}$ was proved by Kollar, Mori, Miyaoka and Takagi [KMMT]. The smooth case was proved by Kollar, Mori and Miyaoka in any dimension. The conjecture is open even in dim 3 when $\delta<1$. In order, in many interesting applications $\delta<1$. The following special case of conjecture \ref{BAB} was proved by Borisov in dim 3 [B] and by McKernan in any dimension [Mc].

\begin{thm}[BM]\label{BM}
The set of all klt WLF pairs $(X,B)$ with a fixed given index is bounded i.e. these pairs are elements of an algebraic family.
\end{thm}

The following conjecture is due to Shokurov.

\begin{conj}[ACC for mlds]\label{acc}
Suppose $\Gamma\subseteq [0,1]$ satisfies the descending chain condition (DCC). Then the following set satisfies the ascending chain condition (ACC):
\[ \{\mld(\mu,X,B)| {\mbox{$(X,B)$ is lc of dim d, $\mu$ a point of $X$ and $B\in \Gamma$}}\}\]

\end{conj}

Note that $\mu$ is assumed to be a Grothendieck point of $X$ which is not necessarily closed. We show the above conjecture as $ACC_{d,\Gamma}$ for short. Alexeev proved $ACC_{2,\Gamma}$ for any DCC set $\Gamma\subseteq [0,1]$ [A2]. This conjecture is open in higher dimensions except in some special cases.

\begin{conj}[Log Termination]\label{lt}
Let $(X,B)$ be a Klt pair of dim d. Then any sequence of $K_X+B$-flips terminates.
\end{conj}

This conjecture is the last step of the LMMP in dim $d$ and we show it as $LT_{d}$. Kawamata proved  $LT_{3}$ [K2] and the four dimensional case with terminal singularities [KMM]. Actually $LT_{4}$ is the main missing component of $LMMP_4$ without which we can not apply the powerful LMMP to problems in algebraic geometry. This conjecture did not seem to be that difficult at least because of the short proof of Kawamata to $LT_{3}$ where he uses the classification of terminal singularities. The later classification is not known in higher dimensions. Recent attempts by Kawamata and others to solve $LT_{4}$ showed that this problem is much deeper than they expected. There is speculation that it may be even more difficult than the flip problem (Shokurov believes this).

We listed several important conjectures with no obvious relation. It is Shokurov's amazing idea to put all these conjectures in a single framework which we call it {\textbf{Shokurov's Program}}: 

\begin{equation}\label{program}
\end{equation}
\begin{description}

\item[${\mathbf{ACC_d \longrightarrow LT_d}}$] Shokurov proved that the $LT_d$ follows from the above ACC conjecture up to dim $d$ and the following problem up to dim $d$ [Sh4]: 

\begin{conj}[Lower Semi-Continuity]\label{lsc}
For any klt pair $(X,B)$ of dim d and any $c\in \{0, 1, \dots, d-1\}$ the function $\mld_c(\mu,X,B): \{c-points ~of~ X\}\longrightarrow \mathbb{R}$ is lower semi-continuous.
\end{conj}

A $c$-point is a $c$-dim Grothendieck point of $X$. This conjecture is proved up to dim 3 by Ambro [Am]. This conjecture doesn't seem to be as tough as previous conjectures. Shokurov proved this problem in dim 4 for mlds in $[0,2]$ [Sh4, lemma 2]. So ACC in dim 4 is enough for the log termination in dim 4 [Sh4, corollary 5]. Actually ACC for mlds in [0,1] for closed points is enough [Sh4, corollary 5].

\item[$\mathbf{BAB_{d-1} \longrightarrow ACC_d}$] Shokurov associates a topological dimension $0\leq reg(P\in X,B)\leq d-1$ to any $d$-dim lc singularity $(P\in X,B)$ [Sh2, 7.9] and proves that the $ACC_{d,\Gamma}$ for pairs with $reg(P\in X,B)=0$ is followed from the BAB conjecture in dim $d-1$ [PSh , 4.4]. In order if $reg(P\in X,B)=0$ then the singularity is exceptional (see definition \ref{exc}). Also $ACC_{d,\Gamma}$ for pairs with $reg(P\in X,B)\in \{1, \dots, d-2\}$ can be reduced to lower dimensions. So the only remaining part of $ACC_{d,\Gamma}$ is when $reg(P\in X,B)=d-1$. This case is expected to be proved using different methods. So {\emph{in particular}} $ACC_{4,\Gamma}$ follows from the BAB in dim 3 and the $reg(P\in X,B)=3$ case. Moreover $ACC_{3,\Gamma}$ is followed from the $reg(P\in X,B)=2$ case. 
 
\item[$\mathbf{WC_{d-1} \longrightarrow BAB_{d-1}}$] And here it comes probably the most important application of the theory of complements: $WC_{\delta, d-1,\{0\}}$ ``implies''  $BAB_{\delta,d-1,[0,1]}$. More precisely, these two problems can be solved together at once. In order, in those situations where boundedness of varieties is difficult to prove, boundedness of complements is easier to prove. And that is exactly what we do in this paper for the 2-dim case: we prove $WC_{\delta, 2,\{0\}}$ and $BAB_{\delta,2,[0,1]}$. Our main objective was to obtain a proof  with as little as possible use of surface geometry so that it can be generalised to higher dimensions. In other words, the methods used in the proof of these results are of the most importance to us. After finishing this work, now, we expect to finish the proof of $WC_{\delta, 3,\{0\}}$ ``implies''  $BAB_{\delta,3,[0,1]}$ in not a far future! And that was our original goal.

\item[The program in dim 4] Let us mention that by carrying out Shokurov's program in dim 4, in which the main ingridient is $WC_{\delta ,3, \{0\}}$ i.e. boundedness of $\epsilon$-lc complements in dim 3, we will prove the following conjectures:
\begin{itemize}
\item ACC for mlds in dim 3.
\item Boundedness of $\delta$-lc 3-fold log Fanos i.e. BAB in dim 3.
\item ACC for mlds in dim 4.
\item Lower semi-continuity for mlds in dim 4.
\item Log termination in dim 4 and then LMMP in dim 4.
\end{itemize}

\end{description}

About this paper:
\begin{enumerate} 
\item Section 1  is devoted to the study of complements on log surfaces.

\item In 1.2 we recall some definitions and lemmas.

\item In 1.3 we prove $WC_{\delta, 1, [0,1]}$ i.e. the boundedness of $\epsilon$-lc complements in dim 1 (\textbf{theorem \ref{curve}}). 

\item In 1.5 we prove $WC_{\delta, 2, \{0\}}$ for the case $X=Z$ i.e. the boundedness of $\epsilon$-lc complements in dim 2, locally, for points on surfaces with $B=0$ (\textbf{theorem \ref{main-local-isom}}). 

\item In 1.6 we prove $WC_{\delta, 2, \{0\}}$ when $X/Z$ is a birational equivalence i.e. the boundedness of $\epsilon$-lc complements in dim 2, locally, for birational contractions of surfaces with $B=0$ (\textbf{theorem \ref{3-I}}). This proof is a surface proof i.e. we heavilly use surface geometry and it is not expected to be generalized to higher dimensions. A second proof of the birational case is given in 1.10 (\textbf{theorem \ref{4-A'}}).

\item In 1.7 we prove $WC_{\delta, 2, \{0\}}$ when $Z=pt$ i.e. the boundedness of $\epsilon$-lc complements on surfaces, globally, with $B=0$ (\textbf{theorem \ref{weak-2dim}}).  The proof  is based on the LMMP and we expect to be generalised to higher dimensions. As a corollary we give a totally new proof to the boundedness of $\epsilon$-lc log Dell Pezzo surfaces (=BAB indim 2) (\textbf{corollary \ref{BAB'}}). Another application of our theorem is that the boundedness of lc ($\epsilon=0$) complements  can be proved only using complement theory (\textbf{theorem \ref{bc_2}}). The later boundedness was proved by Shokurov [Sh2].

\item In 1.8 we give a second proof of $WC_{\delta, 2, \{0\}}$ in the global case i.e. when $Z=pt$ (\textbf{theorem \ref{3-5-H}}). This proof also uses surface geometry and is not expected to be generalized to higher dimensions.

\item In 1.9 we discuss an example which shows that the transformed boundary doesn't  have a better singularity than the original boundary.

\item In 1.10 we give a proof to all local cases for $B\in \Phi_{sm}$, in particular, the case where $X/Z$ is a fibration over a curve (\textbf{theorem \ref{4-A'}}). This proof is also based on the LMMP.

\item Section 2 is about higher dimensional $\epsilon$-lc complements. We discuss our joint work with Shokurov.

\item In 2.1 We give a \textbf{Plan} about how to attack the boundedness of $\epsilon$-lc complements in dimension 3. This is proposed by the author.

\item In 2.2 We give Shokurov's \textbf{Plan} about how to attack the boundedness of $\epsilon$-lc complements in dimension 3. 
 
\end{enumerate}

Let us summarise the main results of section 1:

\begin{thm} Conjecture \ref{weak} holds in dim 1 for $\Gamma=[0,1]$. 
\end{thm}

See \ref{curve} for the proof.

\begin{thm} Conjecture \ref{weak} holds in dim 2 in the global case (i.e. $\dim Z=0$) for $\Gamma=\{0\}$. 
\end{thm}

See \ref{weak-2dim} and \ref{3-5-H} for  proofs.

\begin{thm} Conjecture \ref{weak} holds in dim 2 in the local cases (i.e. $\dim Z>0$) for $\Gamma=\Phi_{sm}$. 
\end{thm}

See \ref{4-A'} , \ref{main-local-isom} and \ref{3-I} for proofs.

\begin{cor} Conjecture \ref{BAB} holds in dim 2.
\end{cor}

See \ref{BAB'} for proof. Conjecture \ref{BAB} in dim 2 was first proved by Alexeev using different methods [A1].

\begin{cor} Theorem \ref{bc_2} can be proved using only the complement theory.
\end{cor}

See the discussion following theorem \ref{bc_2}. 

\begin{rem}[$\epsilon$-lc Complements Method] 
Though formally speaking the list above are the main results in section 1, but we believe that the method used to prove \ref{weak-2dim} and \ref{4-A'} is the most important result of this section. 
\end{rem}

Here we mention some developments in the theory of complements. The following theorem was proved by Shokurov [Sh2] for surfaces.

\begin{thm}\label{bc_2}
 There exists a finite set $\mathcal{N}_{2}$ of positive integers such that any  2-dim lc weak log Fano pair $(X/P\in Z,B)$ has a $(0,n)$-complement$/P\in Z$ for some $n\in \mathcal{N}_{2}$ if $B$ is semi-standard i.e. for each coefficient $b$ of $B$, $b\geq \frac{6}{7}$ or $b=\frac{m-1}{m}$ for some natural number $m$. In order if $\dim Z>0$ then the theorem holds for a general boundary.
\end{thm}

In order Shokurov uses the BAB in dim 2 in the proof of the above theorem. As mentioned before, the results of this paper imply the BAB in dim 2. So  the above  theorem can be proved only based on the theory of complements. A similar theorem is proved by Prokhorov and Shokurov in dim 3 modulo BAB in dim 3 and the effective adjunction in dim 3 [PSh1]. However the local case doesn't need the later assumptions as the following theorem shows.

\begin{thm}
Let $(X/P\in Z,B)$ be a Klt WLF 3-fold pair where $\dim Z\geq 1$ and $B\in \Phi_{sm}$. Then   $K_X+B$ is $(0,n)$-complementary$/P\in Z$ for some $n\in \mathcal{N}_{2}$. 
\end{thm}

Complements have good inductive properties as the theorem above shows which was proved by Prokhorov and Shokurov [PSh]. This theorem  is stated and proved in higher dimensions in more general settings (see [PSh]). To avoid some exotic definitions, we stated only the 3-fold version. 

Finally we give some easy examples of complements. More interesting examples can be found in [Sh1, Sh2, Pr, PSh, PSh1].

\begin{exa}
Let $(X/Z,B)=(\mathbb{P}^1/pt.,0)$ and $P_1, P_2, P_3$ distinct points on  $\mathbb{P}^1$. Then $K_X+P_1+P_2$ is a $(0,1)$-complement for $K_X$ but it is not an $(\epsilon,n)$-complement for any $\epsilon>0$ since $K_X+P_1+P_2$ is not Klt. On the other hand $K_X+\frac{2}{3}P_1+\frac{2}{3}P_2+\frac{2}{3}P_3$ is a $(\frac{1}{3}, 3)$-complement for $K_X$. 
\end{exa}

\begin{exa}
Let $(X_1/Z_1,B_1)=(\mathbb{P}^2/pt.,0)$ and $(X_2/Z_2,B_2)=(\mathbb{P}^2/\mathbb{P}^2,0)$. Then $K_{X_2}$ is a $(2,1)$-complement/$Z_2$ at any point $P\in Z_2$ but obviously $K_{X_1}$ is not even numerically zero/$Z_1$ though $K_{X_1}=K_{X_2}$. 
\end{exa}

\begin{exa}
Let $(X/Z,B)=(X/X,0)$ where $X$ is a surface with canonical singularities. Then the index of $K_X$ is 1 at any point $P\in X$. So we can take $B^+=0$ and $K_X$ a $(1,1)$-complement$/X$ for $K_X$ at any $P\in X$. 
\end{exa}

\subsection*{acknowledgement}

I am deeply grateful to my adviser Professor V.V. Shokurov for introducing me to these problems and helping me through very patiently. I also thanks the EPSRC and the Cecil King Prize for their financial support.

\clearpage

%%%%%%%%%%%%%%%%%%%%%%%%%%%%%%%%%%%%%%%%%%%%%%%%%%%%%%%%%%%%%%%%%%

\subsection{Preliminaries}

In this subsection we bring some basic definitions and constructions.

\begin{defn}
A set $\mathcal{X}$ of varieties of the same dimension is called bounded if there are schemes $\mathbb{X}$ and $S$ of finite type and a morphism $\phi: \mathbb{X}\longrightarrow S$ such that each element of $\mathcal{X}$ is isomorphic to a geometric fibre of $\phi$. Moreover each fibre should give an element of $\mathcal{X}$.
\end{defn}

\begin{defn}
Let $\mathcal{X}$ be a set of pairs $(X,B_X)$ of the same dimension. Then the set consisting of all $(X,\Supp B_X)$ where $X\in \mathcal{X}$ is called bounded if there are schemes $\mathbb{X}$ and $S$ of finite type, a divisor $\mathbb{B}$ on $\mathbb{X}$ and a morphism $\phi: \mathbb{X}\longrightarrow S$ such that each $X\in \mathcal{X}$ is isomorphic to a geometric fibre $\mathbb{X}_{s}$ and $\Supp B_X=\Supp \mathbb{B}|\mathbb{X}_{s}$. And $(X,B)$ bounded means that $(X,\Supp B_X)$ is bounded and there are only finitely many possibilities for the coefficients of $B$. 
\end{defn}

\begin{defn}
Let $(X,B)$ be a Klt pair of dim $d$. Let $\phi: Y\longrightarrow X$ be a morphism such that $B_Y\geq 0$ where $K_Y+B_Y={^*(K_X+B)}$. Then $Y$ is called a partial resolution of $(X,B)$.  
\end{defn}

\begin{lem}
Let $\mathcal{X}=\{X\}$ be a bounded set of Klt varieties of dim $d$ such that $-K_X$ is nef and big. Then $\mathcal{Y}$ the set of partial resolutions for all $X\in\mathcal{X}$ is  bounded. 
\end{lem}
\begin{proof}
By assumptions for each $Y\in\mathcal{Y}$ there are $X\in\mathcal{X}$ and a boundary $B_Y$ such that $K_Y+B_Y={^*K_X}$. Since $\mathcal{X}$ is bounded hence there are only a finite number of possibilities for the coefficients of $B_Y$ independent of $Y$ thus the index of $K_Y+B_Y$ is bounded. And since $-K_X$ is nef and big so $-(K_Y+B_Y)$ is also nef and big. Now by Borisov-Mckernan [Mc] the set $\mathcal{Y}$ is bounded. 
$\Box$
\end{proof}

\begin{defn}
A variety $X/Z$ of dim $d$, is called Pseudo-WLF/$Z$ if there exists a boundary $B$ such that $(X/Z,B)$ is WLF. Moreover $X$ is called Klt Pseudo-WLF/$Z$ if there is a Klt WLF  $(X/Z,B)$. If $\dim Z=0$ then we usually drop $Z$. 
\end{defn}

\begin{rem}
Pseudo-WLF varieties have good properties. For example $\overline{NE}(X/Z)$ is a finite rational polyhedral cone. Moreover each extremal face of the cone is contractible [Sh5, Sh3]. In section 2 we prove that the Klt Pseudo-WLF property is preserved under  flip and divisorial contractions. 
\end{rem}

\begin{defn}
The set $\Phi_{sm}=\{\frac{k-1}{k}| k\in \mathbb{N}\}\cup \{1\}$ is called the set of standard boundary multiplicities. For a boundary $B$ by $B\in \Phi_{sm}$ we mean that the coefficients of $B$ are in $\Phi_{sm}$.
\end{defn}

\begin{defn}[Exceptional pairs]\label{exc}
Let $(X/Z,B)$ be a pair of dim $d$. If $\dim Z=0$ then $(X/Z,B)$ is called exceptional if there is at least a $(0,\mathbb{Q})$-complement $K_X+B^+$ and any $(0,\mathbb{Q})$-complement $K_X+B^+$ is Klt. If $\dim Z>0$ then $(X/Z,B)$ is called exceptional if there is at least a $(0,\mathbb{Q})$-complement $K_X+B^+$ and any $(0,\mathbb{Q})$-complement $K_X+B^+$ is plt on a log terminal resolution. Otherwise $(X/Z,B)$ is called non-exceptional. 
\end{defn}

\begin{rem}\label{analytic-algebraic}
 Boundedness of analytic $(\epsilon,n)$-complements implies the boundedness of  algebraic $(\epsilon,n)$-complement. That is because of the general GAGA priciple [Sh1]. 
\end{rem}

\begin{lem}\label{pre-1}
Let $Y/X/Z$ and $K_Y+B_Y$ be nef$/X$ and $K_X+B={_*(K_Y+B_Y)}$ be $(\epsilon,n)$-complementary/$Z$. Moreover assume that each non-exceptional$/X$ component of $B_Y$ that intersects an exceptional divisor$/X$ has a standard coefficient then $(Y,B_Y)$ will also be $(\epsilon,n)$-complementary/$Z$.
\end{lem}
\begin{proof}
See [PSh, 6.1].
\end{proof}

%%%%%%%%%%%%%%%%%%%%%%%%%%%%%%%%%%%%%%%%%%%%%%%%%%%%%%%%%%%%%%%%

\subsection{{The case of curves}}

In this subsection we prove \ref{weak} for the case of curves. Note that 1-dim global log Fano pairs are just $(\mathbb{P}^{1}, B)$ for a boundary $B=\sum_{i} b_{i}B_{i}$ where $\sum_{i} b_{i}-2<0$. The local case for curves is trivial.  

\begin{thm}\label{curve} $WC_{\delta, 1,[0,1]}$ holds; more precisely, suppose $\frac{m-1}{m}\leq 1-\delta <\frac{m}{m+1}$ for $m$ a natural number then we have:

\begin{description} 
 \item[$\diamond$] $N_{\delta,1, [0,1]}\subseteq \cup_{0<k\leq m}\{k,k+1\}$.
 \item[$\diamond$] $(\mathbb{P}^{1},B^{+})$ can be taken $\frac{1}{m+1}$-lc.
\end{description}
\end{thm}
\begin{proof}
Let $B=\sum_{i} b_{i}B_{i}$ and put $b=b_{h}=max\{b_{i}\}$ and suppose $\frac{k-1}{k}\leq b <\frac{k}{k+1}$ for a natural number $k$. If $k=1$ then $0\leq b <\frac{1}{2}$ and so we have a 1-complement $B^{+}=0$. 

Now assume that $K>1$ and define $a_{i,t}=\llcorner (t+1)b_{i}\lrcorner$ and note that by our assumptions $\sum_{i} a_{i,k}\leq 2k+1$ since $\sum_{i}b_{i}<2$. If $K+B$ doesn't have a $k$-complement then $\sum_{i} a_{i,k}=2k+1$. Since $\frac{k-1}{k}\leq b <\frac{k}{k+1}$ we have $\frac{(k+1)(k-1)}{k}=k+1-\frac{k+1}{k}=k-\frac{1}{k}\leq (k+1)b <\frac{(k+1)k}{k+1}=k$. Thus $a_{h,k}=k-1$ and $1-\frac{1}{k}\leq ~ \langle(k+1)b\rangle ~ <1$ where $\langle .\rangle$ stands for the fractional part. 

Now $a_{i,k+1}=\llcorner (k+2)b_{i}\lrcorner=\llcorner (k+1)b_{i}+b_{i}\lrcorner$. So $a_{i,k+1}$ is $a_{i,k}$ or $a_{i,k}+1$. The later happens iff $1\leq b_{i}+\langle (k+1)b_{i}\rangle$. By the above $b_{h}+\langle (k+1)b_{h}\rangle \geq \frac{k-1}{k}+1-\frac{1}{k} \geq 1$ so  $a_{h,k+1}=a_{h,k}+1$. On the other hand since $\sum_{i}b_{i}<2$ and $\sum_{i}a_{i,k}=2k+1$ then $\sum_{i}\langle (k+1)b_{i}\rangle <1$. And since  $ 1-\frac{1}{k}\leq \langle (k+1)b_{h}\rangle$ then $\langle (k+1)b_{i}\rangle <\frac{1}{k}$ if $i\neq h$. So under assumption $i\neq h$ if $1 \leq \langle (k+1)b_{i}\rangle +b_{i}$ then $1-\frac{1}{k}<b_{i}$. \\
 Thus if $K+B$ has no $k+1$-complement then  $1 \leq \langle (k+1)b_{j}\rangle +b_{j}$ should hold at least for some $j\neq h$. So again we have $ 1-\frac{1}{k}\leq \langle (k+1)b_{j}\rangle$ and so $\langle (k+1)b_{j}\rangle + \langle (k+1)b_{h}\rangle \geq 2(1-\frac{1}{k})\geq 1$ and this is a contradiction. Hence $K+B$ should have a $k$ or $k+1$-complement. If $K+B$ has a $k$ complement then we can have a maximum  $max\{b_{i}^{+}\}=b^{+}=1-\frac{1}{k}\leq 1-\delta$. If it has a $k+1$-complement then again we can have a maximum $b^{+}\leq \frac{k}{k+1}$.  Since $0<k\leq m$ then $N_{\delta,1}\subseteq \cup_{0<k\leq m}\{k,k+1\}$ and $K+B^{+}$ can be chosen as $\frac{1}{m+1}$-lc. These prove the theorem. 
$\Box$
 \end{proof}

\begin{rem} Above we just proved that $\sum_{i} \frac{\llcorner (n+1)b_{i}\lrcorner}{n}\leq 2$ for a bounded $n$. If the equality doesn't hold then  we may add some positive coefficients to get the equality and construct the complement. 
\end{rem}

%%%%%%%%%%%%%%%%%%%%%%%%%%%%%%%%%%%%%%%%%%%%%%%%%%%%%%%%%%%%

\subsection{The case of surfaces}

  We divide the surface case of conjecture \ref{weak} into the following cases: 

\begin{description}
 \item[local isomorphic]  $X/Z$ is an isomorphism.
 \item[local birational]   $X/Z$ is birational but may not be an isomorphism.
 \item[local over curve]  $Z$ is a curve.
 \item[global]   $Z$ is a point.

\end{description}

%%%%%%%%%%%%%%%%%%%%%%%%%%%%%%%%%%%%%%%%%%%%%%%%%%%%%%%%%%%%%

\subsection{Local isomorphic case}

 The main theorem in this subsection is theorem \ref{main-local-isom}. We use classification of surface singularities.

\begin{thm}\label{main-local-isom}
 Conjecture $WC_{\delta, 2,\{0\}}$ holds in the local isomorphic case. 
\end{thm} 

\begin{proof}
 If $\delta\geq 1$ then $P$ is smooth and so we are already done. So assume that $\delta<1$. If the singularity at $P$ is 
of type $E_6$, $E_7$ or $E_8$ then there are only a finite number of possibilities up to analytic isomorphism because of the $\delta$-lc assumption [see Pr 6.1.2]. Otherwise the graph of the resolution will be either of type $A_{r}$:

\begin{displaymath}
\xymatrix{O^{-\alpha_{r}} \ar@{-}[rr]&& \dots &\ar@{-}[rr]&  & &O^{-\alpha_{2  }}\ar@{-}[rr] && O^{-\alpha_{1}}}
\end{displaymath}

 where $\alpha_{i}\geq 2$. Or of type $D_{r}$:  

\begin{displaymath}
\xymatrix{&& &&&& && O^{-2}\ar@{-}[d] \\
 O^{-\alpha_{r}} \ar@{-}[rr]&& \dots &\ar@{-}[rr]& & &O^{-\alpha_{2}}\ar@{-}[rr] && O^{-\alpha_{1}}\\ 
&&&&&&& & O^{-2}\ar@{-}[u]  }
\end{displaymath}

where $\alpha_{i}\geq 2$.

First we work out the case $A_{n}$. Let $K_{W}-\sum_{i}e_{i}E_{i}={^{*}K_{Z}}$ where $e_{i}$ are the discrepancies for a log resolution $W\longrightarrow Z$/$P$. The following lemma is well known and a proof can be found in [AM, 1.2].

\begin{lem} The numbers $(-E_{i}^{2})$ are bounded from above in terms of $\delta$.
\end{lem}

 Intersecting $K_{W}-\sum_{i}e_{i}E_{i}$ with all the exceptional divisors we get a system like the following:

\[a_{1}(-E_{1}^{2})-a_{2}-1= 0\]
\[a_{2}(-E_{2}^{2})-a_{1}-a_{3}= 0   \]
\[a_{3}(-E_{3}^{2})-a_{2}-a_{4}= 0\]
\[\vdots \]
\[a_{r-1}(-E_{r-1}^{2})-a_{r-2}-a_{r}=0 \]
\[a_{r}(-E_{r}^{2})-a_{r-1}-1= 0\]

where $a_i$ is the log discrepancy of $E_i$ with respect to $K_Z$. 

 From the equation $a_{i}(-E_{i}^{2})-a_{i-1}-a_{i+1}\leq 0$ we get 
$a_{i}(-E_{i}^{2}-2)+a_{i}-a_{i-1}\leq a_{i+1}-a_{i}$ which shows that if $ a_{i-1}\leq a_{i}$ then $ a_{i}\leq a_{i+1}$ and moreover if $ a_{i-1}< a_{i}$ then $ a_{i} <a_{i+1}$ . So the solution for the system above should satisfy the following: 

\begin{equation}\label{A_r}
\end{equation}
\[a_{1} \geq \dots \geq a_{i}\leq \dots \leq a_{r}\]

for some $i\geq 1$. If $r\leq 2$ (or any fixed number) then the theorem is trivial. So we may assume that $r>3$ and  also can assume $i\neq r$ unless $a_1=a_2=\dots=a_r$. Now for any $i\leq j<r$ if $-E_{j}^{2}>2$ then $a_{j+1}-a_{j}\geq a_{j}(-E_{j}^{2}-2)\geq \delta$. So if we have $l$ members in $\{j: -E_{j}^{2}>2 ~and ~i\leq j<r\}$ then $a_{r}\geq l\delta$. Hence $a_{r}(-E_{r}^{2}-1)+a_{r}-a_{r-1}\geq l\delta$ which contradicts the last equation in the system if $l$ gets big arbitrarily. In order $l\delta \leq 1$ and so $l \leq \frac{1}{\delta}$. Similar observation shows that the number $l'$ of $1 \leq j\leq i$ where $-E_{j}^{2}>2$ should be bounded. Then $l+l' \leq \frac{2}{\delta}$.

 Suppose $a_{i_{2}}=\dots =a_{i}=\dots =a_{i_{1}}$, $a_{i_{1}-1}\neq a_{i_{1}}$ (or $i_1=1$) and $a_{i_{2}}\neq a_{i_{2}+1}$ (or $i_2=r$) where $i_2\leq i\leq i_1$. Assume that  $i_{1}\neq i$ or  $i_{2}\neq i$. Let for example $i_{1}\neq i$. Then if all $a_j$ are not equal ($=1$) then we have 

\[1=(-E_{r}^{2}-1)a_{r}+a_{r}-a_{r-1}\geq (r-i_{1})(a_{i_{1}+1}-a_{i_{1}})\]
\[=(r-i_{1})[(-E_{i_{1}}^{2}-2)a_{i_{1}}+a_{i_{1}}-a_{i_{1}-1}]\]
\[=(r-i_{1})(-E_{i_{1}}^{2}-2)a_{i_{1}}\geq (r-i_{1})\delta\]
because $-E_{i_{1}}^{2}$ can not be equal to $2$.

So $(r-i_{1})\delta \leq 1$ then $r-i_{1}\leq \frac{1}{\delta}$ is bounded. Similarly $i_{2}$ should be bounded. These observation show that, mentioning that $-E_{k}^{2}$ are bounded, the denominators of $a_{k}$ are bounded. And so the index of $K_{Z}$ at $P$ is bounded and so we are done in this case. But if $i_{1}= i=i_{2}$ then the story is different. In this case note that $\delta \leq (-E_{i}^{2}-2)a_{i}= a_{i-1}-a_{i}+a_{i+1}-a_{i}$. So $\frac{\delta}{2} \leq a_{i-1}-a_{i}$ or $\frac{\delta}{2} \leq a_{i+1}-a_{i}$. For example assume that the later holds then similar to the calculations we just carried out above, $r-i$ will be bounded. But it can happen that $a_{i-1}-a_{i}$ is very small so we won't be able to bound $i$.

 In order, we try to find a solution for the following system with bounded denominators:

\[u_{1}(-E_{1}^{2})-u_{2}-1\leq 0\]
\[u_{2}(-E_{2}^{2})-u_{1}-u_{3}\leq 0   \]
\[u_{3}(-E_{3}^{2})-u_{2}-u_{4}\leq 0\]
\[\vdots \]
\[u_{r-1}(-E_{r-1}^{2})-u_{r-2}-u_{r}\leq 0 \]
\[u_{r}(-E_{r}^{2})-u_{r-1}-1\leq 0\]

To do this, note that if $-E_{i-1}^{2}>2$ then $\delta \leq (-E_{i-1}^{2}-2)a_{i-1}= a_{i-2}-a_{i-1}+a_{i}-a_{i-1}\leq a_{i-2}-a_{i-1}$ then again similar computations to the above shows that $i$ is bounded. Now let $j$ be the smallest number such that $-E_{j}^{2}=\dots = -E_{i-1}^{2}=2$ (remember that we have assumed $\frac{\delta}{2} \leq a_{i+1}-a_{i}$. ). Hence $j$ is bounded. Now take $u_{j}=\dots =u_{i}=\frac{1}{2}$ then the following equations are satisfied if $i-j>2$:

\[u_{j+1}(-E_{j+1}^{2})-u_{j}-u_{j+2}=2u_{j}-u_{j}-u_{j}=0   \]
\[\vdots \]
\[u_{i-1}(-E_{i-1}^{2})-u_{i-2}-u_{i}=2u_{i-1}-u_{i-2}-u_{i}= 0 \]

 Since $r-i$ and $j$ are bounded then the number of remaining equations is bounded and so to satisfy them we just have to divide $u_{i}=\frac{1}{2}$ by a bounded natural number. This finishes the $A_{r}$ type. The $D_{r}$ type will follow.

\begin{rem}
 In order we have constructed a Klt  log divisor $K_W+D$ with bounded index such that  $-(K_W+D)$ is nef and big/$P\in Z$. Now we may use remark \ref{analytic-algebraic}.
\end{rem}

\begin{rem}
All the bounds occurring in the proof are effective and can be calculated in terms of $\delta$.
\end{rem}

\begin{rem}
In Shokurov's case where $\delta=\epsilon=0$ we just take $u_{1}=\dots =u_{r}=0$.
\end{rem}

The case of $D_{r}$: We have a chain $E_{1}, \dots , E_{r}$ of exceptional divisors plus $E$ and $E'$ where $E$ intersects only $E_{1}$ and the same holds for $E'$. In this case we have the following system:

\[a(-E^{2})-a_{1}-1= 0\]
\[a'(-E'^{2})-a_{1}-1= 0\]
\[a_{1}(-E_{1}^{2})-a-a'- a_{2}+1= 0\]
\[a_{2}(-E_{2}^{2})-a_{1}-a_{3}= 0   \]
\[a_{3}(-E_{3}^{2})-a_{2}-a_{4}= 0\]
\[\vdots \]
\[a_{r-1}(-E_{r-1}^{2})-a_{r-2}-a_{r}=0 \]
\[a_{r}(-E_{r}^{2})-a_{r-1}-1= 0\]

Note that $-E^{2}=-E'^{2}=2$ so  $2a-a_{1}-1= 0$ and $2a'-a_{1}-1= 0$ hence $a+a'=a_{1}+1$. Replacing this in the third equation and ignoring the two first equations we get the following system:

\[a_{1}(-E_{1}^{2}-1)-a_{2}= 0\]
\[a_{2}(-E_{2}^{2})-a_{1}-a_{3}= 0   \]
\[a_{3}(-E_{3}^{2})-a_{2}-a_{4}= 0\]
\[\vdots \]
\[a_{r-1}(-E_{r-1}^{2})-a_{r-2}-a_{r}=0 \]
\[a_{r}(-E_{r}^{2})-a_{r-1}-1= 0\]

From this system we get a solution as following:  

 \[a_{1}=\dots =a_{i}<a_{i+1}<\dots <a_{r}\]

(i=r also may happen. In this case $a=a'=a_{1}=\dots =a_{r}=1$).
 Now $r-i$ should be bounded. In order if $i>1$ then $-E_{1}^{2}=\dots =-E_{i-1}^{2}=2$ but $-E_{i}^{2}>2$ (we have assumed $r>i$). Now $\delta (-E_{i}^{2}-2) \leq a_{i}(-E_{i}^{2}-2)+a_{i}-a_{i-1}=a_{i+1}-a_{i}$ (if $i=1$ then $\delta (-E_{1}^{2}-2) \leq a_{1}(-E_{1}^{2}-2)=a_{2}-a_{1}$). We also have the fact that $a_{k+1}-a_{k}\leq a_{k+2}-a_{k+1}$ for $i\leq k <r-1$. And on the other hand $\sum_{i\leq k<r} a_{k+1}-a_{k} \leq a_{r}<a_{r}+a_{r}-a_{r-1}<1$. So we conclude that $r-i$ should be bounded. Moreover since $-E_{k}^{2}$ is bounded, this proves that the denominators of all $a_{k}$ in the $D_{r}$ case are bounded and so the index of $K_{Z}$ at $P$. In this case $B^+=0$ and this finishes the proof of theorem \ref{main-local-isom}. 

$\Box$  
\end{proof}

\begin{rem} Essentially the boundedness properties that we proved and used in the proof of theorem \ref{main-local-isom} had been more or less discovered by some other people independently. Shokurov had used these ideas in an unpublished preprint on mlds. 
\end{rem}

\begin{recall}\label{graphs} Here we recall the diagrams for the $E_{6}$, $E_{7}$ and $E_{8}$ types of singularities. The following is a general case of such singularities:

\begin{displaymath}
\xymatrix{\mathbb{C}^2/{\mathbb{Z}_{m_{1}}} \ar@{-}[r] & O^{-p}& \ar@{-}[l]  \mathbb{C}^2/{\mathbb{Z}_{m_{2}}}\\
 &\ar@{-}[u] O^{-2}& }
\end{displaymath}

where the only possibilities for $(m_{1},m_{2})$ are $(3,3)$, $(3,4)$ and $(3,5)$. So the possible diagrams are as follows: For $(m_{1},m_{2})=(3,3)$ we have

{
\begin{description}

 \item[$1$]   \begin{displaymath}
\xymatrix{ O^{-3} \ar@{-}[r] & O^{-p}& \ar@{-}[l] O^{-3}\\
 &\ar@{-}[u] O^{-2}& }
\end{displaymath}

 \item[$2$]   \begin{displaymath}
\xymatrix{O^{-2} \ar@{-}[r] & O^{-2} \ar@{-}[r] & O^{-p}& \ar@{-}[l] O^{-3}\\
 &&\ar@{-}[u] O^{-2}& }
\end{displaymath}
 
\item[$3$]   \begin{displaymath}
\xymatrix{O^{-2} \ar@{-}[r] & O^{-2} \ar@{-}[r] & O^{-p}& \ar@{-}[l] O^{-2}& \ar@{-}[l] O^{-2}\\
 &&\ar@{-}[u] O^{-2}&& }
\end{displaymath}

\end{description}
}

For $(m_{1},m_{2})=(3,4)$ we have

\begin{description}

\item[$4$]    \begin{displaymath}
\xymatrix{ O^{-3} \ar@{-}[r] & O^{-p}& \ar@{-}[l] O^{-4}\\
 &\ar@{-}[u] O^{-2}& }
\end{displaymath} 

\item[$5$]   \begin{displaymath}
\xymatrix{O^{-2} \ar@{-}[r] & O^{-2} \ar@{-}[r] & O^{-p}& \ar@{-}[l] O^{-4}\\
 &&\ar@{-}[u] O^{-2}& }
\end{displaymath}

\item[$6$]   \begin{displaymath}
\xymatrix{ O^{-3} \ar@{-}[r] & O^{-p}& \ar@{-}[l] O^{-2}& \ar@{-}[l] O^{-2}& \ar@{-}[l] O^{-2}\\
 &\ar@{-}[u] O^{-2}&&& }
\end{displaymath}

\item[$7$]   \begin{displaymath}
\xymatrix{O^{-2} \ar@{-}[r] & O^{-2} \ar@{-}[r] & O^{-p}& \ar@{-}[l] O^{-2}& \ar@{-}[l] O^{-2}& \ar@{-}[l] O^{-2}\\
 &&\ar@{-}[u] O^{-2}&&& }
\end{displaymath}

\end{description}

And for $(m_{1},m_{2})=(3,5)$ we have  

{
\begin{description}

\item[$8$]    \begin{displaymath}
\xymatrix{ O^{-3} \ar@{-}[r] & O^{-p}& \ar@{-}[l] O^{-5}\\
 &\ar@{-}[u] O^{-2}& }
\end{displaymath} 

\item[$9$]   \begin{displaymath}
\xymatrix{O^{-2} \ar@{-}[r] & O^{-2} \ar@{-}[r] & O^{-p}& \ar@{-}[l] O^{-5}\\
 &&\ar@{-}[u] O^{-2}& }
\end{displaymath}

\item[$10$]   \begin{displaymath}
\xymatrix{ O^{-3} \ar@{-}[r] & O^{-p}&  \ar@{-}[l] O^{-2}& \ar@{-}[l] O^{-3}\\
 &\ar@{-}[u] O^{-2}&& }
\end{displaymath}

\item[$11$]   \begin{displaymath}
\xymatrix{O^{-2} \ar@{-}[r] & O^{-2} \ar@{-}[r] & O^{-p}& \ar@{-}[l] O^{-2}& \ar@{-}[l] O^{-3} \\
 &&\ar@{-}[u] O^{-2}&& }
\end{displaymath}

\item[$12$]   \begin{displaymath}
\xymatrix{O^{-3}  \ar@{-}[r] & O^{-p}& \ar@{-}[l] O^{-3}& \ar@{-}[l] O^{-2} \\
 &\ar@{-}[u] O^{-2}&& }
\end{displaymath}

\item[$13$]   \begin{displaymath}
\xymatrix{O^{-2}  \ar@{-}[r] &O^{-2}  \ar@{-}[r]& O^{-p}& \ar@{-}[l] O^{-3}& \ar@{-}[l] O^{-2} \\
 &&\ar@{-}[u] O^{-2}&& }
\end{displaymath}

\item[$14$]   \begin{displaymath}
\xymatrix{O^{-3}  \ar@{-}[r] & O^{-p}& \ar@{-}[l] O^{-2}& \ar@{-}[l] O^{-2}& \ar@{-}[l] O^{-2} & \ar@{-}[l] O^{-2}\\
 &\ar@{-}[u] O^{-2}&& &&}
\end{displaymath}

\item[$15$]   \begin{displaymath}
\xymatrix{O^{-2}  \ar@{-}[r] &O^{-2}  \ar@{-}[r]& O^{-p}& \ar@{-}[l] O^{-2}& \ar@{-}[l] O^{-2}& \ar@{-}[l] O^{-2}& \ar@{-}[l] O^{-2} \\
 &&\ar@{-}[u] O^{-2}&&&& }
\end{displaymath}

\end{description}
}

\end{recall}

\clearpage

%%%%%%%%%%%%%%%%%%%%%%%%%%%%%%%%%%%%%%%%%%%%%%%%%%%%%%%%%%%%%

\subsection{Local birational case}

In this subsection wherever we write /$Z$ we mean /$P\in Z$ for a fixed point $P$ on $Z$.

\begin{thm}\label{3-I} $WC_{\delta, 2,\{0\}}$ holds in the birational case.
\end{thm}

{\textbf{Strategy of the proof:}} Let $W$ be a minimal resolution of $X$ and $\{E_{i}\}$, $\{F_{j}\}$ be the exceptional divisors /$Z$ on $W$ where the $E_{i}$ are exceptional/$X$ but $ F_{j}$ are not ($E$ will be used for a typical $E_{i}$ and similarly $F$ for $F_{j}$ or its birational transform). We will construct an antinef/$Z$ and Klt log divisor $K_{W}+\Omega=K_{W}+\sum_{i} u_{i}E_{i}+\sum_{j} u_{j}F_{j}$  where $u_{i},u_{j}<1$ are rational numbers with bounded denominators. Then we use remark \ref{analytic-algebraic}.

\begin{proof} 

 By contracting those curves where $-K_{X}$ is numerically zero, we can assume that $-K_{X}$ is ample/$Z$ (we can pull back the complement). Let $W$ be the minimal resolution of  $X$. Then since $K_{W}$ is nef/$X$ by the negativity lemma we have $K_{W}-\sum_{i}e_{i}E_{i}=K_{W}+\sum_{i}(1-a_{i})E_{i} \equiv {^*K_{X}}$ where $ e_{i}\leq 0$.

\begin{defn} For any smooth model $Y$ where $W/Y/Z$ we define $\overline{exc}(Y/Z)$ to be the graph of the exceptional curves ignoring the birational transform of exceptional divisors of type $F$. For an exceptional/$Z$ divisor $G$ on $Y$ not of type $F$, $\overline{exc}(Y/Z)_{G}$ means the connected component of $\overline{exc}(Y/Z)$ where $G$ belongs to. 
\end{defn}

\begin{lem} We have the followings on $W$:
\begin{description}
 \item[$\diamond$] The exceptional divisors/$Z$ on $W$ are with simple normal crossings.
 \item[$\diamond$] Each $F$ (i.e. each exceptional divisor of type $F$) is a $-1$-curve.
 \item[$\diamond$] The model $\overline{W}$ obtained by blowing down $-1$-curves/$Z$ is the minimal resolution of $Z$.
 \item[$\diamond$] Each $F$ cuts at most two exceptional divisors of type $E$. 
\end{description} 
\end{lem}

\begin{proof}
 Let $F$ be an exceptional divisor/$Z$ on $W$ which is not exceptional/$X$. Then $(K_{W}-\sum_{i}e_{i}E_{i}).F= K_{W}.F+\sum_{i}(-e_{i})E_{i}.F=2p_{a}(F)-2-F^{2}+\sum_{i}(-e_{i})E_{i}.F < 0$ where $p_{a}(F)$ stands for the arithmetic genus of the curve $F$. Then $2p_{a}(F)-2-F^{2} < 0$ and so $p_{a}(F)=0$ and $-F^{2}=1$. In other words $F$ is a $-1$-curve.\\
On the other hand by contracting $-1$-curves/$Z$ (i.e. running the classical minimal model theory for smooth surfaces on $W/Z$) we get a model $\overline{W}/Z$ where $K_{\overline{W}}$ is nef/$Z$. Actually $\overline{W}$ is the minimal resolution of $P\in Z$. 

The exceptional divisors/$Z$ on $\overline{W}$ are with simple normal crossings and since $W$ is obtained from $\overline{W}$ by a sequence of blow ups then the exceptional divisors/$Z$ on $W$ also would be with simple normal crossings. This observation gives more information. Since all the $F$, exceptional/$Z$ but not/$X$, are contracted/$\overline{W}$ then they should intersect at most two of $E_{i}$ because $exc(\overline{W}/Z)$ is with simple normal crossings and $F$ is exceptional/$\overline{W}$. Moreover no two exceptional divisors of type $F$ should intersect on $W$ because they are both $-1$-curves. This means that the intersection points of any two exceptional divisor/$Z$ on $X$ should be a singular point of $X$. Also any exceptional divisor/$Z$ on $X$ contains at most two singular points of $X$. 
$\Box$

\end{proof}

Now Suppose $\{Q_{k}\}_{k}$ to be the singular points of $X$. If no one of the points $\{Q_{k}\}$ is of type $A_{r}$ then the proof of theorem \ref{main-local-isom} shows that the discrepancies $e_{i}$ are with bounded denominators so we are already done. But if there is one point of type $A_{r}$ then the proof is more complicated (surprisingly the $A_{r}$ type is the most simple case in the sense of Shokurov i.e. when $\delta=0$). Similar to the proof of theorem \ref{main-local-isom} we try to understand the structure of $exc(W/Z)$ and the blow ups $W\rightarrow \overline{W}$.

\begin{defn} A smooth model $\ddot{W}$ where $W/\ddot{W}$ and $\ddot{W}/\overline{W}$ are series of smooth blow ups, is called a $blow~ up~ model$ of $\overline{W}$. Such a model is called $minimal$ if there is $X'$ such that $K_{\ddot{W}}$ is nef/$X'$ and $X/X'/Z$. In other words it is the minimal resolution of $X'$. The connected components of $\overline{exc}(\ddot{W}/Z)$ are either of type $A_{r}$, $D_{r}$, $E_{6}$, $E_{7}$ or $E_{8}$ for a minimal blow up model.
\end{defn}

\begin{defn}\label{3-N}  We call the divisor $K_{W}+\omega=K_{W}+\sum_{i} (1-a_{i})E_{i}={^*K_{X}}$  the $primary$ $~log~ divisor$. The pair $(X,B)$ has a  log canonical $n$-complement $K_{X}+B^+$ over $Z$ in the sense of Shokurov [Sh2] ($n\leq 6$). From now on we call it a $Shokurov~ complement$. So  $K_{W}+\omega_{Sh}+C=K_{W}+\sum_{i}(1-a^{Sh}_{i})E_{i}+\sum_{j}(1-a^{Sh}_{j})F_{j}+C={^*(K_{X}+B^+)}$ where $C$ is the birational transform of the non-exceptional part of $B^+$. We call  $K_{W}+\omega_{Sh}$ a $Shokurov~ log~ divisor$ and the numbers $a^{Sh}_{i}$ and $a^{Sh}_{j}$ Shokurov log discrepancies.   
\end{defn}

\begin{defn}\label{3-F} In the graph $exc(W/Z)$ if we ignore those $F$ which appear with zero coefficient in $\omega_{Sh}$ (i.e. $a^{Sh}=1$) then we get a graph $exc(W/Z)_{>0}$ with some connected components. The connected graph $\mathcal{C}$ consisting of exceptional/$Z$ curves with $a^{Sh}=0$, is in one of the components of the graph $exc(W/Z)_{>0}$ which we show by $\mathcal{G}$ ($\mathcal{C}$ is connected because of the connectedness of the locus of log canonical centres/$P\in Z$). Now contracting all $-1$-curves/$Z$ in $\mathcal{G}$ and continuing the contractions of subsequent $-1$-curves/$Z$ which appear in $\mathcal{G}$, finally we get a model which we show as $W_{\mathcal{G}}$. The transform of $\mathcal{G}$ on $W_{\mathcal{G}}$ is shown by $\mathcal{G}_{1}$ and similarly the transform of $\mathcal{C}$ is $\mathcal{C}_{1}$.
\end{defn}

\begin{defn}\label{3-H} A chain of exceptional curves  consisting of $G_{\beta_{1}}, \dots, G_{\beta_{r}}$ is called $strictly ~monotonic$ if $r=1$ or if $a_{\beta_{1}}< a_{\beta_{2}}<\dots < a_{\beta_{r}}$ (these are log discrepancies with respect to $K_X$). $G_{\beta_{1}}$ is called the $base~ curve$.
\end{defn}

\begin{defn} Let $G\in exc(\ddot{W}/Z)$ for a smooth blow up model $\ddot{W}$. Then define the $negativity$ of $G$ on this model as $N_{\ddot{W}}(G)=(K_{\ddot{W}}+{_{*}\omega}).G\leq 0$ ($_{*}\omega$ is pushdown of $\omega$). We also define the $total~ negativity$ by $N_{\ddot{W}}=\sum_{\alpha} N_{\ddot{W}}(G_{\alpha})$ where $G_{\alpha}$ runs on all exceptional divisors/$Z$ on $\ddot{W}$ and for $G\in \overline{exc}(\ddot{W}/Z)$ define $N_{\ddot{W},G}=\sum_{\alpha} N_{\ddot{W}}(G_{\alpha})$ when the sum runs on all members of $\overline{exc}(\ddot{W}/Z)_{G}$. Similarly define the negativity functions $N^{Sh}$ and $N^+$ replacing $\omega$ with $\omega_{Sh}$ and $\omega_{Sh}+C$ respectively. Note that the later is always zero, because $K_{W}+\omega_{Sh}+C\equiv 0/Z$.     

\end{defn}

\begin{defn} The (smooth) blow up of a point which belongs to two exceptional divisors/$Z$ on a model is called a $double~ blow~ up$. If this point belongs to just one exceptional divisor/$Z$ then we call it a $single~ blow~ up$. A blow up is called  double$^+$ blow up if the blown up point belongs to two components of $_{*}(\omega_{Sh}+C)$ (this is the pushdown). Similarly define a $single^+$ blow up.

\end{defn}

\begin{lem}\label{3-A}  For any exceptional $G_{\beta}\in \overline{exc}(\ddot{W}/Z)$ on a blow up model $\ddot{W}$ we have:
\begin{description}
 \item[$\diamond$] $-1+\delta \leq N_{\ddot{W},G_{\beta}}$ if $\overline{exc}(\ddot{W}/Z)_{G_{\beta}}$ is of type $D_{r}$, $E_{6}$, $E_{7}$ or $E_{8}$. In particular, in these cases  $-1+\delta \leq N_{\ddot{W}}(G_{\beta})$ holds.
 \item[$\diamond$] $2(-1+\delta) \leq N_{\ddot{W},G_{\beta}}$ and $-1+\delta \leq N_{\ddot{W}}(G_{\beta})$ if $\overline{exc}(\ddot{W}/Z)_{G_{\beta}}$ is of type $A_{r}$ except if it is strictly monotonic.   
\end{description}
\end{lem}
\begin{proof} $D_{r}$ case: Similar to the notation in the proof of theorem \ref{main-local-isom} let $G_{\beta}, G_{\beta'}, G_{\beta_{1}}, \dots , G_{\beta_{r}}$ be the exceptional divisors in $\overline{exc}(\ddot{W}/Z)_{G_{\beta}}$. Then from the equations in the proof of theorem \ref{main-local-isom} for the $D_{r}$ case we get the following system for the log discrepancies:

\[2a_{\beta}-a_{\beta_{1}}-1\leq 0\]
\[2a_{\beta'}-a_{\beta_{1}}-1\leq 0\]
\[2a_{\beta_{1}}-a_{\beta}-a_{\beta'}- a_{\beta_{2}}+1\leq 0\]
\[2a_{\beta_{2}}-a_{\beta_{1}}-a_{\beta_{3}}\leq 0   \]
\[\vdots \]
\[2a_{\beta_{r-1}}-a_{\beta_{r-2}}-a_{\beta_{r}}\leq 0 \]
\[2a_{\beta_{r}}-a_{\beta_{r-1}}-1\leq 0\]

Adding the first and the second equations gives $2a_{\beta}+2a_{\beta'}-2a_{\beta_{1}}-2\leq 0$ and putting this in the third equation we get $a_{\beta_{1}}\leq a_{\beta_{2}}$ and so  $~a_{\beta_{1}}\leq a_{\beta_{2}}\leq \dots \leq a_{\beta_{r}}$. So
  
 \[N_{\ddot{W},G_{\beta}}\geq a_{\beta}+a_{\beta'}+a_{\beta_{r}}-a_{\beta_{1}}-2\geq a_{\beta}+a_{\beta'}+a_{\beta_{2}}-a_{\beta_{1}}-2\]\[ \geq 2a_{\beta_{1}}+1-a_{\beta_{1}}-2\geq a_{\beta_{1}}-1\geq \delta -1\] because $2a_{\beta_{1}}+1\leq a_{\beta}+a_{\beta'}+a_{\beta_{2}}$ and the fact that $X$ is $\delta$-lc.

 $A_{r}$ case [non-strictly monotonic]: In this case assume that the exceptional divisors in $\overline{exc}(\ddot{W}/Z)_{G_{\beta}}$ are $ G_{\beta_{1}}, \dots , G_{\beta_{r}}$ so we get the system: 

\[2a_{\beta_{1}}-a_{\beta_{2}}-1\leq 0\]
\[2a_{\beta_{2}}-a_{\beta_{1}}-a_{\beta_{3}}\leq 0   \]
\[\vdots \]
\[2a_{\beta_{r-1}}-a_{\beta_{r-2}}-a_{\beta_{r}}\leq 0 \]
\[2a_{\beta_{r}}-a_{\beta_{r-1}}-1\leq 0\]

 So there will be $k$ such that $a_{\beta_{1}}\geq a_{\beta_{2}}\geq \dots \geq a_{\beta_{k}}\leq a_{\beta_{r}}$. Thus $N_{\ddot{W}}(G_{\beta_{1}})\geq a_{\beta_{1}}+a_{\beta_{1}}-a_{\beta_{2}}-1\geq  a_{\beta_{1}}-1\geq \delta-1$. In this way we get the similar inequalities for all other equations  except for $N_{\ddot{W}}(G_{\beta_{k}})$. Suppose  $N_{\ddot{W}}(G_{\beta_{k}})<\delta-1$. So we get  $2a_{\beta_{k}}-a_{\beta_{k-1}}-a_{\beta_{k+1}}< \delta-1$ and so $1- \delta <a_{\beta_{k-1}}+a_{\beta_{k+1}}-2a_{\beta_{k}}\leq a_{\beta_{1}}+a_{\beta_{r}}-2a_{\beta_{k}}$. On the other hand by adding all the equations in the system above we get $N_{\ddot{W},G_{\beta}}\geq a_{\beta_{1}}+a_{\beta_{r}}-2>1- \delta+2a_{\beta_{k}}-2\geq \delta-1$. This is a contradiction with the fact that  $N_{\ddot{W}}(G_{\beta_{k}})\geq N_{\ddot{W},G_{\beta_{k}}}$. 

To get the inequality for $ N_{\ddot{W},G_{\beta_{k}}}$ add all the equations in the system above. Note that if $r=2$ then $a_{\beta_{1}}=a_{\beta_{2}}$ and checking the lemma is easy in this case.

$E_{6}$, $E_{7}$, $E_{8}$ {cases:\footnote{I just prove that $-1+\delta \leq N_{\ddot{W}}(G)$ for any exceptional $G$. We won't need the inequality for total negativity.}}  In these cases the graph   $\overline{exc}(\ddot{W}/Z)_{G_{\beta}}$ is as in the recall \ref{graphs}. It is enough to put $2$ in place of all self-intersection numbers because the negativity becomes smaller. We start from the smallest possible graph i.e. Case 1 in the mentioned recall:

\[2a_{\beta}-a_{\beta_{2}}-1\leq 0\]
\[2a_{\beta_{2}}-a_{\beta}-a_{\beta_{1}}- a_{\beta_{3}}+1\leq 0\]
\[2a_{\beta_{1}}-a_{\beta_{2}}-1\leq 0\]
\[2a_{\beta_{3}}-a_{\beta_{2}}-1\leq 0\]

Adding all the equations we get $N_{\ddot{W},G_{\beta}}=a_{\beta}+a_{\beta_{1}}+a_{\beta_{3}}-a_{\beta_{2}}-2$. And by the second equation we have $a_{\beta}+a_{\beta_{1}}+a_{\beta_{3}}-a_{\beta_{2}}\geq a_{\beta_{2}}+1$ so  $N_{\ddot{W},G_{\beta}}\geq a_{\beta_{2}}+1-2\geq \delta-1$. In order this was a special case of the $D_{r}$ type (the similarity of the system not necessarily the graph $\overline{exc}(\ddot{W}/Z)_{G_{\beta}}$). Note that the inequality for the total negativity implies the inequality for the negativity of each exceptional curve.\\
 Now we prove other cases by induction on the number of the exceptional curves. The minimum is $4$ exceptional curves and we just proved this case. Suppose we have proved up to $k-1$ and our graph has $k$ members. Suppose the exceptional curves are $G_{\beta}, G_{\beta_{1}}, \dots ,G_{\beta_{k-1}}$ and such that  $G_{\beta_{l}}$ cuts $G_{\beta}, G_{\beta_{l-1}}$ and $G_{\beta_{l+1}}$. If $l=2$ or $l=k-2$ then again this will be a system of type $D_{r}$. Assume that otherwise happens. Since $ -a_{\beta}+1\geq 0$ so we get a system as follows 

\[2a_{\beta_{1}}-a_{\beta_{2}}-1\leq 0\]
\[2a_{\beta_{2}}-a_{\beta_{3}}-a_{\beta_{1}}\leq 0\]
\[\vdots\]
\[2a_{\beta_{k-1}}-a_{\beta_{k-2}}-1\leq 0\]

This is a system of type $A_{k-1}$ so either we have $a_{\beta_{1}}\geq a_{\beta_{2}}$
 or $a_{\beta_{k-1}}\geq a_{\beta_{k-2}}$. Assume that the first one holds (the other case is similar). Now note that $N_{\ddot{W}}(G_{\beta_{1}})\geq 2a_{\beta_{1}}-a_{\beta_{2}}-1=a_{\beta_{1}}-a_{\beta_{2}}+a_{\beta_{1}}-1\geq \delta-1$. By ignoring $G_{\beta_{1}}$ we get a system for a graph with a smaller number of elements:

\[2a_{\beta_{2}}-a_{\beta_{3}}-1\leq 2a_{\beta_{2}}-a_{\beta_{3}}-a_{\beta_{1}} \leq 0\]
\[\vdots\]
\[2a_{\beta_{l}}-a_{\beta_{l-1}}-a_{\beta_{l+1}}-a_{\beta}+1\leq 0\]
\[\vdots\]
\[2a_{\beta_{k-1}}-a_{\beta_{k-2}}-1\leq 0\]

 and the lemma is proved by induction. 
$\Box$
\end{proof}

\begin{lem}\label{3-B} Suppose $\xi \in \ddot{W}/\overline{W}$ ($\ddot{W}$ is a blow up model), ~$\tilde{W}$ the blow up of $ \ddot{W}$ at $\xi$ and $G_{\alpha}$ the exceptional divisor of the blow up. Then we have the followings: 

If $G_{\alpha}$ is the double blow up of $G_{\beta}$ and $G_{\gamma}$ (i.e. $\xi \in G_{\beta}\cap G_{\gamma}$) then:

\begin{description} 
 \item[$\diamond$] $N_{\tilde{W}}(G_{\alpha})=a_{\alpha}-a_{\beta}-a_{\gamma}$ where $a_{\alpha}$ is the log discrepancy of $G_{\alpha}$ for $K_{X}$ and similarly $a_{\beta}$ and $a_{\gamma}$. 
 \item[$\diamond$] $N_{\tilde{W}}(G_{\beta})=N_{\ddot{W}}(G_{\beta})-N_{\tilde{W}}(G_{\alpha})$ and  $N_{\tilde{W}}(G_{\gamma})=N_{\ddot{W}}(G_{\gamma})-N_{\tilde{W}}(G_{\alpha})$. 
 \item[$\diamond$] $N_{\tilde{W}}=N_{\ddot{W}}-N_{\tilde{W}}(G_{\alpha})$.
\end{description}

And if $G_{\alpha}$ is the single blow up of $G_{\beta}$ then 
\begin{description} 
 \item[$\diamond$]  $N_{\tilde{W}}(G_{\beta})=N_{\ddot{W}}(G_{\beta})-N_{\tilde{W}}(G_{\alpha})$, $N_{\tilde{W}}(G_{\alpha})=a_{\alpha}-a_{\beta}-1\leq -\delta$ and $N_{\ddot{W}}(G_{\beta})+\delta\leq 0$.
 \item[$\diamond$] $N_{\tilde{W}}=N_{\ddot{W}}$.
\end{description}

\end{lem}
\begin{proof}
Standard computations and left to the reader.
$\Box$
\end{proof}

\begin{cor}\label{3-C} If $G_{\alpha}$ is a single blow up of $G_{\beta}$ on $\ddot{W}$, a blow up model of $\overline{W}$, and  $N_{\ddot{W}}(G_{\beta})\geq \delta -1$ then $a_{\alpha}\geq a_{\beta}+\delta$. 
\end{cor}
\begin{proof} Since $G_{\alpha}$ is a single blow up of $G_{\beta}$ then $1+a_{\beta}-a_{\alpha}+N_{\ddot{W}}(G_{\beta})\leq 0$ and so $1+a_{\beta}-a_{\alpha}+\delta -1 \leq 0$ then $a_{\beta}+\delta \leq a_{\alpha}$.
\end{proof}

\begin{defn} Let $\xi$ be a point on a blow up model $\ddot{W}$. Define the \\ $multiplicity~  of~ double~  blow~  ups$ as 

\[\mu_{db}(\xi)=max \{\#\{~ double~ blow~ ups~/\xi~ before~ having~ a~single~ blow ~up~/\xi\}\}\] 

the maximum is taken over all sequences of blow ups from $\ddot{W}$ to $W$. The next lemma shows the boundedness of this number. 
\end{defn}

\begin{lem}\label{3-E} $\mu_{db}(\xi)$ is bounded. 
\end{lem}
\begin{proof} Since by lemma \ref{3-B} each double blow up adds a non-negative number to the total negativity of the system and since the total negativity is {bounded\footnote{because the total negativity on $\overline{W}$ is bounded. This boundedness for the $A_{r}$ and $D_{r}$ cases is shown in lemma \ref{3-A} and for other cases it is obvious}}  then except a bounded number of double blow ups we have $\frac{-\delta}{2}\leq N_{\tilde{W}}(G_{\alpha})=a_{\alpha}-a_{\beta}-a_{\gamma}\leq 0$ where $G_{\alpha}$ is the double blow up of some $G_{\beta}$ and $G_{\gamma}$ and  $G_{\beta}\cap G_{\gamma}/\xi$. The inequality shows that $a_{\beta}+\frac{\delta}{2}\leq a_{\beta}+a_{\gamma}-\frac{\delta}{2}\leq a_{\alpha}$ and similarly $a_{\gamma}+\frac{\delta}{2}\leq a_{\alpha}$. In other words the log discrepancy is increasing at least by $\frac{\delta}{2}$. Since log discrepancies are in $[\delta,1]$ then the number of these double blow ups has to be bounded. 
$\Box$
\end{proof}

\begin{defn} Let $\xi \in \ddot{W}$ a blow up model. Define the  $single~ blow~ up  \\ ~multiciplity$ ~of ~$\xi$ as:

 \[\mu_{sb}(\xi)=max\{\#\{G: G ~is ~single~ blown~ up ~and ~ G/\xi\}\}\] 

The maximum is taken over all sequences of blow ups from $\ddot{W}$ to $W$. In the above definition $G$ is the excpetional divisor of a sinlge blow up/$\xi$. Also define $ \mu_{sb}(G_{\beta})=\sum_{\xi\in G_{\beta}} \mu_{sb}(\xi)$ and $ \mu_{sb}(\ddot{W})=\sum_{\xi\in \ddot{W}} \mu_{sb}(\xi)$ .
\end{defn}

So if $\xi_{2}/\xi_{1}$ (these points may be on different models) then $\mu_{sb}(\xi_{1})\geq \mu_{sb}(\xi_{2})$.

\begin{rem} Usually there is not a unique sequence of blow ups from $\ddot{W}$ to $W$. In order if $\xi_{1}\neq \xi_{2}$ are points on $\ddot{W}$ and they are centres of some exceptional divisors on $W$ then it doesn't matter which one we first blow up to get to $W$ i.e. they are independent.   
\end{rem}

\begin{defn} Let $\xi \in \overline{exc}(\ddot{W}/Z)$ be a point on a blow up model $\ddot{W}$. We call such a point a $base~ point$ if there is an exceptional divisor $G_{\alpha}/\xi$ on a blow up model $\tilde{W}$  such that $N_{\tilde{W}}(G_{\alpha})<\delta -1$.
\end{defn}

\begin{rem}\label{3-K}  By lemma \ref{3-A} and lemma \ref{3-B} if $\xi \in \overline{exc}(\ddot{W}/Z)_{G_{\beta}}$ and $\overline{exc}(\ddot{W}/Z)_{G_{\beta}}$ is of type $A_{r}$ (non-strictly monotonic) , $D_{r}$, $E_{6}$, $E_{7}$ or $E_{8}$ then $\xi$ can not be a base point. Moreover again by lemma \ref{3-A} if $\overline{exc}(\ddot{W}/Z)_{G_{\beta}}$ is strictly monotonic then there can be at most one base point in $\overline{exc}(\ddot{W}/Z)_{G_{\beta}}$ and it only can belong to the base curve. 
\end{rem}

\begin{lem}\label{3-M} $\mu_{sb}(\xi)$ is bounded if $\xi\in \ddot{W}$  is not a base point.  
\end{lem}

\begin{proof}  If $ G_{\alpha}/\xi$ is a single blown up exceptional divisor then since $\xi$ is not a base point we get $\delta-1\leq N_{\tilde{W}}(G_{\alpha})$. So if $G_{\alpha})$ is a single blow up/$\xi$ of $G_{\beta})$ then $a_{\alpha}\geq a_{\beta}+\delta$ i.e. it  increases the log discrepancy at least by $\delta$. And as we said in the proof of lemma \ref{3-E}, except a bounded number, any other double blow up/$\xi$ increases the log discrepancy at least by $\frac{\delta}{2}$. In order, there can be only a bounded number of blow ups/$\xi$ from $\ddot{W}$ to $W$.
$\Box$

\end{proof}

\begin{cor}\label{3-L}  The number of exceptional curves/$\xi$ on $W$ is bounded for any non-base point $\xi\in \ddot{W}$. 
\end{cor}

We now continue to the proof of theorem \ref{3-I}. If no divisor in $\omega_{Sh}$ has coefficient 1 then this is what we are looking for. Since in this case $K_{W}+\omega_{Sh}$ will be a $1/6$-lc log divisor. If the opposite happens i.e. some divisor appear with coefficient 1  in $\omega_{Sh}$ then these divisors will form a connected chain $\mathcal{C}$ which doesn't intersect any other exceptional divisor/$Z$ with positive coefficient in $\omega_{Sh}$, except in the edges of this chain. Some of the exceptional divisors of type $F$ may appear with positive coefficients and some with zero coefficients in $\omega_{Sh}$.

 The image of the graph $\mathcal{G}$ on $W_{\mathcal{G}}$, that is $\mathcal{G}_{1}$ (see definition \ref{3-F}), is either of type $A_{r}$, $D_{r}$, $E_{6}$, $E_{7}$ or $E_{8}$ because similar to what we proved above for $\overline{W}$ the model $W_{\mathcal{G}}$ is the minimal resolution of some surface i.e. the minimal resolution of the surface $X_{\mathcal{G}}$ obtained from $X$ by contracting the exceptional/$Z$ curves  on $X$ whose birational trasnform belong to $\mathcal{G}$. In order there is no $-1$-curve/$X_{\mathcal{G}}$ on $W_{\mathcal{G}}$.

 Now suppose $\mathcal{G}_{1}$ is of type $A_{r}$, not strictly monotonic and let the pushdown of the chain $\mathcal{C}$ be $\mathcal{C}_{1}$ on $W_{\mathcal{G}}$ . Let the exceptional divisors of $\mathcal{G}_{1}$ be $G_{\beta_{1}}, \dots , G_{\beta_{r}}$ and assume that the chain $\mathcal{C}_{1}$ consists of $G_{\beta_{k}}, \dots , G_{\beta_{l}}$. Hence $N^{Sh}_{ W_{\mathcal{G}}}(G_{\beta_{k}})\leq -\frac{1}{6}$,  $N^{Sh}_{ W_{\mathcal{G}}}(G_{\beta_{k+1}})= \dots = N^{Sh}_{ W_{\mathcal{G}}}(G_{\beta_{l-1}})=0$ and  $N^{Sh}_{ W_{\mathcal{G}}}(G_{\beta_{l}})\leq -\frac{1}{6}$ . Here the superscript $Sh$ means that we compute the negativity according to the Shokurov log divisor not the primary log divisor (we have already defined it above). Note that if $a^{Sh}_{\beta}>0$ for some $\beta$ then $a^{Sh}_{\beta}\geq \frac{1}{6}$ because the denominator of $a^{Sh}_{\beta}$ is in $\{1,2,3,4,6\}$. The chain $\mathcal{C}_{1}$ is of type $A_{l-k+1}$. From the constructions in the local isomorphism subsection we can replace  the Shokurov log numbers $a^{Sh}_{\beta_{k}}=0,\dots , a^{Sh}_{\beta_{l}}=0$ with new log numbers with bounded denominators and preserve all other Shokurov log numbers in the graph $exc(W_{\mathcal{G}}/Z)$ such that we obtain a new log divisor $K_{W_{\mathcal{G}}}+\Omega_{1}$ on $W_{\mathcal{G}}$ which is antinef/$Z$ and is a Klt log divisor. Now put $K_{W}+\Omega=^*(K_{W_{\mathcal{G}}}+\Omega_{1})$. The only problem with $\Omega$ is that it may have negative coefficients (it is a subboundary). Remark \ref{3-K} and corollary \ref{3-L} assure us that the negativity of these coefficients is bounded from below. Moreover if an exceptional divisor has negative coefficient in $\Omega$ then it should belong to the graph $\mathcal{G}$. But any exceptional divisor in  $\mathcal{G}$ appears with positive coefficient in $\omega_{Sh}$. Since $\omega_{Sh}\geq \omega$ and also by the definition of  $\mathcal{G}$, any exceptional divisor of type $F$   in $\mathcal{G}$ has positive coefficient at least $\frac{1}{6}$. And if $E$ is not of type $F$ but belongs to $\mathcal{G}$ then since $B^+$ is not zero $P\in Z$ we get positive coefficients in $\omega_{Sh}$ for all exceptional/$Z$ curves which are not of type $F$. Thus all members of $\mathcal{G}=Exc(W/Q)$ appear with positive coefficient in $\omega_{Sh}$. 

Now consider the sum 

\[K_{W}+\Omega+I[K_{W}+\omega_{Sh}]=(1+I)K_{W}+[\Omega+I\omega_{Sh}]\] 

where $I$ is an integer. Mentioning the fact that the negative coefficients appeared in $\Omega$ are bounded from below, this implies that there is a large bounded $I$ such that the sum $\Omega+I\omega_{Sh}$ is an effective divisor. So by construction the log divisor  $K_{W}+\frac{[\Omega+I\omega_{Sh}]}{1+I}$ is $\epsilon$-lc and antinef/$Z$ for some fixed rational number $0<\epsilon$ and the denominators of the coefficients in the log divisor are bounded.

 Now assume that $\mathcal{G}_{1}$ is strictly monotonic and the base curve is $G_{\beta_{1}}$.  By corollary \ref{3-L} and remark \ref{3-K} the only place where we may have difficulties is the base point, $\xi$, on the base curve if there is any such point.

 Now  we blow $\xi$  up and get the exceptional divisor $G_{\alpha_{1}}$. The chain $G_{\alpha_{1}}, G_{\beta_{1}}, \dots , G_{\beta_{r}}$ is not exactly of type $A_{r+1}$ because $G_{\alpha_{1}}$ is a $-1$-curve. But still we can claim that there is at most a base on this chain and it only can be on $G_{\alpha_{1}}$. Obviously a base point cannot be on $G_{\beta_{2}}, \dots , G_{\beta_{r}}$. Now suppose that the intersection point of $G_{\alpha_{1}}$ and $G_{\beta_{r}}$ is a base point. Then the sum of negativities of all $G_{\alpha_{1}}, G_{\beta_{1}}, \dots , G_{\beta_{r}}$ should be less than $2\delta-2$. This is impossible because the sum of negativities of all $G_{\beta_{1}}, \dots , G_{\beta_{r}}$ on $W_{\mathcal{G}}$ is at least $2\delta-2$ (remember that blowing up reduces negativity). 

Now if on $G_{\alpha_{1}}$ there is a base point $\xi_{1}$ then again we blow this point up to get  $G_{\alpha_{2}}$ and so on. This process has to stop after finitely many steps (not after bounded steps!).  Let the final model to be $W_{\xi}$ and $G_{\alpha_{1}}, \dots , G_{\alpha_{s}}$ the new exceptional divisors. In order we have constructed a  chain (because on each curve there was at most one base point) and by adding the new exceptional divisors to $\mathcal{G}_{1}$ we get a new graph $\mathcal{G}_{2}$. Now there is no base point on $\mathcal{G}_{2}$. All the divisors  $G_{\alpha_{i}}$ have self-intersection equal to $-2$ except $G_{\alpha_{s}}$ which is a $-1$-curve. 

Now let $\mathcal{C}_{2}$ to be the pushdown of $\mathcal{C}$ i.e. the connected chain of curves with coefficient one in $\omega_{Sh}$ on $W_{xi}$. If $G_{\alpha_{s}}$ is not in $\mathcal{C}_{2}$ then we proceed exactly as in the non-monotonic case above; that is we assign appropriate coefficients to the members of $\mathcal{C}_{2}$ and keep all other coefficients in $\omega_{Sh}$ on $W_{xi}$. If $G_{\alpha_{s}}$ is in $\mathcal{C}_{2}$ then let $\mathcal{C}'$ be the chain $\mathcal{C}_{2}$ except the member $G_{\alpha_{s}}$. This new chain (i.e. $\mathcal{C}'$) if of type $A_x$ and so we can assign appropriate coefficients to its members and put the coefficient of $G_{\alpha_{s}}$ simply equal to zero and keep all other coefficients in $\omega_{Sh}$ on $W_{xi}$. In any case we construct a Klt log divisor $K+\Omega$ on $W_{xi}$ which is antinef/$Z$ and the boundary coefficients are with bounded denominators. The rest is as in the non-monotonic case above.

Suppose the graph $\mathcal{G}_{1}$ is of type $D_{r}$ and $\mathcal{C}{1}\neq \emptyset$ (if it is empty then we already have $\Omega_{1}$). Assume that the members of $\mathcal{G}_{1}$ are $G_{\beta}, G_{\beta'}, G_{\beta_{1}}, \dots , G_{\beta_{r}}$ and the members of  $\mathcal{C}_{1}$ are  $G_{\beta_{k}}, \dots , G_{\beta_{l}}$. As in the proof of lemma \ref{3-A} for the $D_{r}$ case we have $a^{Sh}_{\beta_{1}}\leq a^{Sh}_{\beta_{2}}\leq \dots$ . So $k=1$ and  we have $2a^{Sh}_{\beta}-0-1\leq 0$ and so $a^{Sh}_{\beta}\leq \frac{1}{2}$ and similarly  $a^{Sh}_{\beta'}\leq \frac{1}{2}$. The chain $\mathcal{C}_{1}$ is of type $A_l$ and so we can change the coefficients of its members in $\omega_{Sh}$ on $W_{\mathcal{G}}$. The rest of the argument is very similar to the above cases. Just note that there is no base point in this case.

 The cases  $E_{6}$, $E_{7}$ and $E_{8}$ are safe by remark \ref{3-K} and corollary \ref{3-L}. In these cases the graph $\mathcal{G}$ is bounded so assigning the primary log numbers to the members of $\mathcal{G}_{1}$ and Shokurov log numbers to the rest of the graph $exc(W_{\mathcal{G}}/Z)$ gives a log divisor which can be used as $K_{W_{\mathcal{G}}}+\Omega_{1}$. Here the proof of theorem \ref{3-I} is finished. 

$\Box$
 
\end{proof}

%%%%%%%%%%%%%%%%%%%%%%%%%%%%%%%%%%%%%%%%%%%%%%%%%%%%%%%%%%%%%%%%%%%%

\subsection{Global case}

The main theorem of this subsection is the following theorem. A generalised version of this and the BAB follow as  corollaries. 

\begin{thm}\label{weak-2dim}
 Conjecture $WC_{\delta, 2,\{0\}}$ holds in the global case i.e. when $Z$ is a point.

\end{thm}

\begin{proof}
We divide the problem into two main cases: exceptional and non-exceptional. $(X,0)$ is non-exceptional if there is a non-Klt $\mathbb{Q}$-complement $K_X+M$. By [Sh2, 2.3.1], under our assumptions on $X$, non-exceptionallity is equivalent to the fact that $K_X$ has a non-Klt $(0,n)$-complement for some $n<58$. We prove that the exceptional cases are bounded. But in the non-exceptional case we only prove the existence of an $(\epsilon,n)$-complement for a bounded $n$. Later we show that this in order implies the boundedness of $X$.

First assume that $(X,0)$ is {\textbf{non-exceptional}}.

\begin{enumerate}

\item Lets show the set of accumulation points of the mlds in dim 2 for lc pairs $(T,B)$ where $B\in \Phi_{sm}$, by $Accum_{2,\Phi_{sm}}$. Then  $Accum_{2,\Phi_{sm}}\cap [0,1]=\{1-z\}_{z\in \Phi_{sm}}=\{\frac{1}{k}\}_{k\in\mathbb{N}}\cup \{0\}$ [Sh8]. Now if there is a $\tau>0$ such that $\mld(P,T,B)\notin [\frac{1}{k},\frac{1}{k}+\tau]$ for any natural number $k$ and any point $P\in T$ then there will be only  a finite number of possibilities for the index of $K_{T}+B$ at $P$ if $(T,B)$ is $\frac{1}{m}$-lc for some $m\in \mathbb{N}$. Now Borisov-Mckernan [Mc, 1.2] implies the boundedness of all such $T$ if $-(K_T+B)$ is nef and big and $\tau$ and $m$ are fixed. In order in the following steps we try to reduce our problem to this situation in some cases. 

\item    
 
\begin{defn}\label{D-tau}
Let $B=\sum b_{i}B_{i}$ be a boundary on a variety $T$ and $\tau>0$ a real number. Define 
\[D_{\tau}:=\sum_{b_{i}\notin [\frac{k-1}{k}-\tau, \frac{k-1}{k}]}b_{i}B_{i}+\sum_{b_{i}\in [\frac{k-1}{k}-\tau, \frac{k-1}{k}]} \frac{k-1}{k}B_{i}\] 

where in the first term  $b_{i}\notin [\frac{k-1}{k}-\tau, \frac{k-1}{k}]$ for any natural number $k$ but in the second term $k$ is the smallest natural number satisfying $b_{i}\in [\frac{k-1}{k}-\tau, \frac{k-1}{k}]$.
\end{defn}

\begin{lem}\label{transform} For any natural number $m$ there is a real number $\tau>0$ such that if $(T,B)$ is a surface log pair, $P\in T$, $K_T+B$  is $\frac{1}{m}$-lc at $P$ and $D_{\tau}\in\Phi_{sm}$  then  $K_T+D_{\tau}$ is also $\frac{1}{m}$-lc at $P$. 
\end{lem}

Note that $\tau$ depends only on $m$. 

\begin{proof}
 By applying the ACC to all surface pairs with standard boundary, we get a fixed rational number $v>0$ such that if any $K_T+D_{\tau}$ is  not $\frac{1}{m}$-lc at $P$ then $\mld(P,T,D_{\tau})<\frac{1}{m}-v$. 

Now assume that the lemma is not true. So there is a sequence $\tau_{1}>\tau_{2}>\dots$ and a sequence of pairs $\{(T_i,B_i)\}$ where if we take $\tau_i$ for the pair $(T_i,B_i)$ then the lemma doesn't hold at $P_i\in T_i$. In other words  $\mld(P_i,T_i,D_{\tau_i})<\frac{1}{m}-v$. 

Write $B_i:=F_i+C_i$ where $F_i=\sum f_{i,x}F_{i,x}$ and $C_i=\sum c_{i,y}C_{i,y}$ have no common components and the coefficient of any component of $C_i$ is equal to the coefficient of the same component in $D_{\tau_i}$ but the coefficient of any component of $F_i$ is less than the coefficient of the same component in $D_{\tau_i}$.

 Now there is a set $\{s_{1,x}\}\subseteq [\frac{m-1}{m}-\tau_{1}, \frac{m-1}{m}]$ of rational numbers such that $\mld(P_1,T_1, \sum s_{1,x}F_{1,x}+C_1)=\frac{1}{m}-v$. There is $i_2$ such that $\max \{s_{1,x}\}<\frac{m-1}{m}-\tau_{i_2}$. So there is also a set $\{s_{2,x}\}\subseteq [\frac{m-1}{m}-\tau_{i_2}, \frac{m-1}{m}]$ such that $\mld(P_{i_2},T_{i_2}, \sum s_{2,x}F_{i_2,x}+C_{i_2})=\frac{1}{m}-\frac{v}{2}$. By continuing this process we find $\{s_{j,x}\}\subseteq [\frac{m-1}{m}-\tau_{i_j}, \frac{m-1}{m}]$ such that $\max \{s_{i_{j-1},x}\}<\frac{m-1}{m}-\tau_{i_j}$. Hence we can find a set $\{s_{j,x}\}\subseteq [\frac{m-1}{m}-\tau_{i_j}, \frac{m-1}{m}]$ such that $\mld(P_{i_j},T_{i_j}, \sum s_{j,x}F_{i_j,x}+C_{i_j})=\frac{1}{m}-\frac{v}{j}$. 

In order we have constructed a set $\cup \{s_{j,x}\}$ of rational numbers which satisfies the DCC condition but there is an increasing set of mlds corresponding to boundaries with coefficients in $\cup \{s_{j,x}\}$. This is a contradiction with the ACC for mlds. 
$\Box$
\end{proof}

 \item  Let $m$ be the smallest number such that $\frac{1}{m}\leq \delta$. Let $h=\min\{\frac{k-1}{k}-\frac{u}{r!}>0\}_{1\leq k\leq m} $ where $u,k$ are natural numbers and $r=\max\{m,57\}$. Now choose a $\tau$ for $m$ as in lemma \ref{transform} such that $\tau <h$. 

Blow up one exceptional divisor $E$ via $f:Y\longrightarrow X$ such that the log discrepancy satisfies $\frac{1}{k}\leq a(E,X,0)\leq \frac{1}{k}+\tau$ for some $k>1$ (if such $E$ doesn't exist then go to step 1). The crepant log divisor $K_Y+B_Y$ is $\frac{1}{m}$-lc and so by lemma \ref{transform} $K_Y+D_{\tau}$ is also $\frac{1}{m}$-lc ($D_{\tau}$ is constructed for $B_Y$). Let $K_X+B^+$ be a $(0,n)$-complement for some $n<58$ and $K_Y+B_{Y}^+$ be the crepant blow up. Then by the way we chose $\tau$ we have $D_{\tau}\leq B^+$. Now run the anti-LMMP over $K_Y+D_{\tau}$ i.e. contract any birational type extremal ray $R$ such that $(K_Y+D_{\tau}).R>0$. At the end of this process we get a model $X_1$ and the corresponding map $g:Y\longrightarrow X_1$. After contracting those birational extremal rays where $K_{X_1}+D_{\tau}$ is numerically zero we get a model  $S_1$ with one of the following properties:

\begin{description}
  \item[$\diamond$] $\rho(S_1)=1$ and $K_{S_1}+D_{\tau}\equiv K_{S_1}+B_{S_1}^+\equiv 0$ and $\frac{1}{m}$-lc.
  \item[$\diamond$]  $\rho(S_1)=2$ and $(K_{S_1}+D_{\tau}).R= 0$ for a non-birational type extremal ray $R$ on $S_1$ and $K_{S_1}+D_{\tau}$ is $\frac{1}{m}$-lc.
  \item[$\diamond$] $-(K_{S_1}+D_{\tau})$ is nef and big and $K_{S_1}+D_{\tau}$ is $\frac{1}{m}$-lc.
\end{description}

where $K_{S_1}+D_{\tau}$ is the birational transform of $K_{Y}+D_{\tau}$.

In any case $-(K_{S_1}+D_{\tau})$ is nef because $D_{\tau}\leq B_{S_1}^+$ and so $D_{\tau}$ can not be positive on a non-birational extremal ray. $K_{S_1}+D_{\tau}$ is $\frac{1}{m}$-lc by the way we have chosen $\tau$.

\item If the first case occurs in the division in step 3 then we are done.

\item If the second case occurs in the division in step 3 then $R$ defines a fiberation $\phi: S_1\longrightarrow Z$. Note that $B_{S_1}^+=D_{\tau}+N$ where each component of $N$ is a fibre of $\phi$ and there are only a finite number of possibilities for the coefficients of $N$. Now we can replace $N$ by $N'\equiv N$ where each component of $N'$ is a general fibre of $\phi$, with only a finite number of possibilities for the coefficients of $N'$ and such that $K_{S_1}+D_{\tau}+N'$ is $\frac{1}{m}$-lc. Note that the components of $N'$ are smooth curves and intersect the components of $D_{\tau}$ transversally in smooth points of $S_1$. Now the only problem is that we don't know if the index of $K_{S_1}+D_{\tau}+N'$ is bounded or not. Note that it is enough if we can get the boundedness of the index of $K_{S_1}+D_{\tau}$.

\item Now assume that the third case or the second case occurs in the division in step 3. Let $C$ be a curve contracted by $g:Y\longrightarrow X_1$ constructed in step 3. If $C$ is not a component of $B_Y$ then the log discrepancy of $C$ with respect to $K_{X_1}+B_{X_1}$ is at least 1 where $K_{X_1}+B_{X_1}$ is the birational transform of $K_{Y}+B_{Y}$. Moreover $g(C)\in \Supp B_{X_1}\neq \emptyset$. So the log discrepancy of $C$ with respect to $K_{X_1}$ is more than 1. This means that $C$ is not a divisor on a minimal resolution $W_1\longrightarrow X_1$. Let $W\longrightarrow X$ be a minimal resolution. Then there is a morphism $W\longrightarrow W_1$. Hence $exc(W_1/X_1)\subseteq exc(W/X)$. Now if $C\in exc(W/X)$ is exceptional/$X_1$ then $a(C,X_1,D_{\tau})<a(C,X,0)$.      

\item Let $(X_1,B_1):=(X_1,D_{\tau})$ and repeat the process. In other words again we blow up one exceptional divisor $E$ via $f_1:Y_1\longrightarrow X_1$ such that the log discrepancy satisfies $\frac{1}{k}\leq a(E,X_1,B_1)\leq \frac{1}{k}+\tau$ for some natural number $k>1$. The crepant log divisor $K_{Y_1}+B_{1,Y_1}$ is $\frac{1}{m}$-lc and so by lemma \ref{transform} $K_{Y_1}+D_{1,\tau}$ is $\frac{1}{m}$-lc. Note that the point which is blown up on $X_1$ can not be smooth since $\tau <h$ as defined in step 3. So according to step 6 the blown up divisor $E$ is a member of $exc(W/X)$. Now we again run the anti-LMMP on $K_{Y_1}+D_{1,\tau}$ and proceed as in step 3.

$$  \xymatrix{
    W \ar[d]\ar[r] & W_1 \ar[d]\ar[r] & W_2 \ar[d]\ar[r] & \dots \\
    Y \ar[d]^{f}\ar[rd]^{g}  &  Y_1\ar[d]^{f_1}\ar[rd]^{g_1} & Y_2 \ar[d]\ar[rd] &\dots\\
    X & X_1\ar[d] & X_2\ar[d] & \dots \\
     & S_1 & S_2 & \dots \\
}$$

\item Steps 6,7 show that each time we blow up a member of $exc(W/X)$ say $E$. And if we blow that divisor down in some step then the log discrepancy $a(E,X_j,B_j)$  will decrease. That divisor will not be blown up again unless the log discrepancy drops at least by $\frac{1}{2(m-1)}-\frac{1}{2m}$. So after finitely many steps either the case one occurs in the division in step 3 or we get a model $X_i$ with a standard boundary $B_i$ such that  there is no $E$ where $\frac{1}{k}\leq a(E,X_i,B_i)\leq \frac{1}{k}+\tau$ for any $1<k\leq m$. The later implies the boundedness of the index of $K_{X_i}+B_i=K_{X_i}+D_{i-1,\tau}$. If $-(K_{X_i}+B_{i})$ is nef and big (case one) then $(X_i,B_i)$ will be  bounded by step 1. Otherwise we have the second case in the division above and so by step 5 we are done (the index of $K_{X_i}+D_{i-1,\tau}+N'$ is bounded).

Now we treat the {\textbf{exceptional}} case: From now on we assume that $(X,0)$ is exceptional.

 \item Let $W\longrightarrow X$  be a minimal resolution. Let $0<\tau<\frac{1}{2}$ be a  number and the minimal log discrepancy of $(X,0)$ be $a=\mld(X,0)$. If $a\geq \frac{1}{2}+\tau$ then we know that $X$ belongs to a bounded family  according to step 1 above. So we assume $a<\frac{1}{2}+\tau$ and then blow up an exceptional/$X$ curve $E_1$ with log discrepancy $a_{E_1}=a(E_1,X,0)\leq \frac{1}{2}+\tau$ to get $Y\longrightarrow X$ and put $K_{Y}+B_{Y}={^*K_{X}}$. Let $t\geq 0$ be a number such that there is an extremal ray $R$ such that $(K_{Y}+B_{Y}+tE_1).R=0$ and $E_1.R>0$ ( and s.t. $K_{Y}+B_{Y}+tE_1$ Klt and antinef). Such $R$ exists otherwise there is a $t>0$ such that $K_{Y}+B_{Y}+tE_1$ is lc (and not Klt) and antinef. This is a contradiction  by [Sh2, 2.3.1]. Now contract $R: Y\longrightarrow Y_1$ if it is of birational type.

Again by increasing $t$ there will be an extremal ray $R_1$ on $Y_1$ such that $(K_{Y_1}+B_{Y_1}+tE_1).R_1=0$ and $E_1.R_1>0$ (preserving the nefness of $-(K_{Y_1}+B_{Y_1}+tE_1)$ ). If it is of birational then contract it and so on. After finitely many steps we get a model $(V_1, B_{V_1}+t_1E_1)$ and a number $t_1>0$ with the following possible outcomes:

\begin{equation}\label{3-5-J}
\end{equation}
\begin{description}
 \item[$\diamond$] $ (V_1, B_{V_1}+t_1 E_1)$ is Klt, $\rho(V_1)=1$ and $K_{V_1}+B_{V_1}+t_1 E_1$ is antinef.
 \item[$\diamond$] $ (V_1, B_{V_1}+t_1 E_1)$ is Klt and $\rho(V_1)=2$ and there is a non-birational extremal ray $R$ on $V_1$. Moreover $K_{V_1}+B_{V_1}+t_1 E_1$ and $K_{V_1}$ are antinef.
 \item[$\diamond$] $ (V_1, B_{V_1}+t_1 E_1)$ is Klt and $\rho(V_1)=2$ and there is a non-birational extremal ray $R$ on $V_1$. Moreover $K_{V_1}+B_{V_1}+t_1 E_1$ is antinef but $K_{V_1}$ is not antinef. 
\end{description}

Define $K_{V_1}+D_1=K_{V_1}+B_{V_1}+t_1 E_1$. Note that in all the cases above $E_1$ is a divisor on $V_1$ and the coefficients of $B_{V_1}$ and $D_1$ are $\geq \frac{1}{2}-\tau$.

\begin{lem}\label{bound-index}
Let $P\in U$ be a $\delta$-lc surface singularity. Moreover suppose that there is at most one exceptional/$U$ divisor such that $a(E, U,0)<\frac{1}{2}+\tau$. Then the index $K_U$ is bounded at $P$ where the bound only depends on $\delta$ and $\tau$. 
\end{lem}

\begin{proof}
We only need to prove this when the singularity is of type $A_r$ (otherwise the index is bounded). If there is no $E/P$ such that $a(E,U,0))<\frac{1}{2}+\frac{\tau}{2}$ then step 1 shows that the index is bounded. But if there is one $E/P$ such that $a(E,U,0)<\frac{1}{2}+\frac{\tau}{2}$ then using the notation as in \ref{A_r} we have $a_{i+1}-a_i\geq \frac{\tau}{2}$ and $a_{i-1}-a_i\geq \frac{\tau}{2}$. This implies the boundedness of $r$ and so the index at $P$.
$\Box$
\end{proof}

\item  Let $U$ be a surface with the following properties:

\begin{description}
 \item[$\diamond$] $\rho(U)=1$.
 \item[$\diamond$] $K_U+G_U$  antinef, Klt and exceptional.
 \item[$\diamond$] $K_U$ antiample.
\end{description}

Now blow up two divisors $E$ and $E'$ as $f: Y_U\longrightarrow U$ such that $a(E,U,0)<\frac{1}{2}+\tau$ and $a(E',U,0)<\frac{1}{2}+\tau$ (suppose there are such divisors). Choose $t,t'\geq 0$ such that $(f^*(K_{U}+G_U)+tE+t'E').R=0$ for an extremal ray $R$ s.t. $R.E\geq 0$ and $R.E'\geq 0$ and $f^*(K_{U}+G_U)+tE+t'E'$ is antinef and Klt. We contract $R$ to get $g:Y_U\longrightarrow U'$. We call such operation a {\textbf{hat of first type}}. Note that $E$ and $E'$ are divisors on $U'$ and $\rho(U')=2$. Define $K_{U'}+G_{U'}$ to be the pushdown of $f^*(K_{U}+G_U)+tE+t'E'$. 

If $K_U$ is $\delta$-lc and such $E,E'$ don't exist as above then the index of $K_U$ will be bounded by lemma \ref{bound-index}. So $U$ will be bounded.

\item  Let $U$ be a surface with the following properties:

\begin{description}
 \item[$\diamond$] $\rho(U)=2$.
 \item[$\diamond$] $K_U+G_U$  antinef, Klt and exceptional.
 \item[$\diamond$] $-K_U$ nef and big.
\end{description}

Now blow up a divisor $E$ to get $f:Y_U\longrightarrow U$ such that $a(E,U,0)<\frac{1}{2}+\tau$ (suppose there is such $E$). Let $t\geq 0$ be such that $(f^*(K_{U}+G_U)+tE).R=0$ for an extremal ray $R$ s.t. $R.E\geq 0$ and $f^*(K_{U}+G_U)+tE$ is antinef and Klt. We contract $R$  to get $g:Y_U\longrightarrow U'$. We call such operation a {\textbf{hat of second type}}.  Note that $E$ is a divisor on $U'$ and $\rho(U')=2$. Define $K_{U'}+G_{U'}$ to be the pushdown of $f^*(K_{U}+G_U)+tE$.

If $K_U$ is $\delta$-lc and such $E$ doesn't exist as above then the index of $K_U$ and so $U$ will be bounded by lemma \ref{bound-index}.

\item  Let $U$ be a surface with the following properties:

\begin{description}
 \item[$\diamond$] $\rho(U)=2$ and $U$ is Pseudo-WLF.
 \item[$\diamond$] There is a birational type extremal ray $R_{bir}$ and the other extremal ray of $U$ is of fibration type.
 \item[$\diamond$] $K_U+G_U$  antinef, Klt and exceptional.
 \item[$\diamond$] $K_U.R_{bir}>0$.
\end{description}

 Then we say that $U$ is of {\textbf{2-bir}} type.  Let $C$ be the divisor that defines $R_{bir}$ on $U$. There is a $c\in (0,1)$ such that $(K_U+cC).C=0$. Now blow up $E$ as $Y_U\longrightarrow U$ such that $a(E,U,cC)<\frac{1}{2}+\tau$ (suppose there is such $E$). Now let $t\geq 0$ such that $f^*(K_{U}+G_U+tC).R=0$ for an extremal ray $R$ s.t. $ R.E\geq 0$, $R.C\geq 0$ and $f^*(K_{U}+G_U+tC)$ is antinef and Klt. We contract $R$ to get $g:Y_U\longrightarrow U'$. We call such operation a {\textbf{hat of third type}}. Define $K_{U'}+G_{U'}$ to be the pushdown of $f^*(K_{U}+G_U+tC)$. Note that in this case $E$ and $C$ are both divisors on $U'$ and $\rho(U')=2$.

If $K_U+cC$ is $\delta$-lc and such $E$ doesn't exist as above then contract $C: U\longrightarrow U_1$. Thus the index of $K_{U_1}$ will be bounded at each point by lemma \ref{bound-index} and so $U_1$ and consequently $U$ will be bounded.

$$  
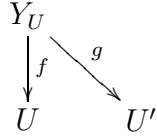
\xymatrix{
    Y_U \ar[d]^{f}\ar[rd]^{g}  &  \\
    U & U' \\
}$$

\item Let $U$ be a surface such that $\rho(U)=2$ and $K_U+G_U$ antinef, Klt and exceptional where $G_U\neq 0$. Moreover suppose there are two exceptional curves $H_1$ and $H_2$ on $U$. In this case let $C$ be a component of $G_U$ and let $t\geq 0$ such that $(K_U+G_U+tC).H_i=0$ for $i=1$ or 2 and $K_U+G_U+tC$  Klt and antinef (assume $i=1$). We contract $H_1$ as $U\longrightarrow U_1$ and define $K_{U_1}+G_{U_1}$ to be the pushdown of $K_{U}+G_U+tC$.

\begin{defn} 
Define $K_U+\Delta_U$ as follows: $K_U+\Delta_U:=K_U$ in step 10 and step 11. $K_U+\Delta_U:=K_U+cC$ in step 12. And $K_{U_1}+\Delta_{U_1}:=K_{U_1}$ in step 13.
\end{defn}

\item The following lemmas are crucial to our proof.

\begin{lem}\label{bound-1}
Let $\mathcal{U}$ be a bounded family of surfaces with Picard number one or two and let $0<x<1$ be a rational number. Moreover assume the following for each member $U_i$:

\begin{description}
 \item[$\diamond$] $-(K_{U_i}+B_i)$ is nef and big for a boundary $B_i$ where each coefficient of $B_i$ is $\geq x$.  
 \item[$\diamond$] $K_{U_i}+B_i$ is Klt.
\end{description}
 
Then $(U_i, \Supp B_i)$ is bounded.

\end{lem}
\begin{proof}
 In order we prove that there is a finite set $\Lambda_{f}$ such that for each $U_i$ there is a boundary $M_i\in \Lambda_f$ s.t. $-(K_{U_i}+M_i)$ is nef and big and $M_i\leq B_i$.

If $\rho(U_i)=1$ then simply take $M_i= x\sum_{\alpha}B_{\alpha}$ where $B_i=\sum_{\alpha}B_{\alpha}$. Obviously $-(K_{U_i}+M_i)$ is nef and big and since $U_i$ belongs to a bounded family so $(U_i, \Supp M_i)$ is also bounded. 

Now suppose $\rho(U_i)=2$. Put $N_i= x\sum_{\alpha}B_{\alpha}$. If $-(K_{U_i}+N_i)$ is not nef then there should be an exceptional curve $E$ on $U_i$ where $(K_{U_i}+N_i).E>0$. Let $\theta: U_i\longrightarrow U'_i$ be the contraction of $E$. By our assumptions $K_{U'_i}+B'_i$, the pushdown of  $K_{U_i}+B_i$, is antiample. So $K_{U'_i}+N'_i$, the pushdown of $K_{U_i}+N_i$ is also antiample. Boundedness of $U_i$ implies the boundedness of $U'_i$ (since we have a bound for the Picard number of a minimal resolution of $U'_i$). Thus $-(K_{U_i}+M_i):=-\theta^*(K_{U'_i}+N'_i)=-(K_{U_i}+N_i+yE)$ is nef and big and there are only a finite number of possibilities for $y>0$. This proves the boundedness of $(U_i, \Supp(N_i+yE))$. Note that in the arguments above $\Supp B_i=\Supp M_i$. 

$\Box$
\end{proof}

\begin{lem}[The main lemma]\label{mainlemma}
 Suppose that  $\mathcal{U}=\{(U,\Supp D)\}$ is a bounded family of log pairs of dim $d$ where $K_U+D$ is antinef and $\epsilon$-lc for a fixed $\epsilon>0$. Then the set of partial resolutions of all $(U,D)\in\mathcal{U}$ is a bounded family.
\end{lem}

Note that here we don't assume $(U,D)$ to be bounded i.e. the coefficients of $D$ may not be in a finite set.

\begin{proof}

 Let $(U_t,D_t)$ be a member of the family. By our assumptions the number of components of $D_t$ is bounded (independent of $t$) and so we can consider any divisor supported in $D_t$ as a point in a real finite dimensional space. Let $D_t=\sum_{1\leq i \leq q} d_{i,t}D_{i,t}$ and define 

 \[\mathcal{H}_t:=\{(h_1,\dots ,h_q)\in \mathbb{R}^q ~|~K_{U_t}+\sum_{1\leq i \leq q} h_{i}D_{i,t}~is~ antinef ~and ~ \epsilon-lc  \}\]

So $\mathcal{H}_t$ is a subset of the cube $[0,1]^q$ and since being $\epsilon$-lc and antinef are closed conditions then $\mathcal{H}_t$ is a closed and hence compact subset of $[0,1]^q$. In oder $\{(U_t,\mathcal{H}_t)\}$ is a bounded family. For each $H\in \mathcal{H}_t$ the corresponding pair $(U_t, H)$ is $\epsilon$-lc. Let $Y_H\longrightarrow U_t$ be a terminal blow up of $(U_t, H)$ and assume that the set of exceptional/$U_t$ divisors on $Y_H$ is $R_H$. For different $H$ we may have different $R_H$ but the union of all $R_H$ is a finite set where $H$ runs through $\mathcal{H}_t$. Suppose otherwise  so there is a sequence $\{H_1,\dots, H_m,\dots\}\subseteq\mathcal{H}_t $ such that the union of all $R_{H_i}$ is not finite. Since $\mathcal{H}_t$ is compact then there is at least an accumulation point in $\mathcal{H}_t$, say $\bar{H}$, for the sequence (we can assume that this is the only accumulation point). So $(U_t, \bar{H})$ is $\epsilon$-lc. Let $v=(1,\dots,1)\in \mathbb{R}^q$. Then there are $\alpha, \beta>0$ such that $K_{U_t}+H_\alpha$ is $\epsilon-\beta$-lc where $\epsilon-\beta >0$ and $H_\alpha$ is the corresponding divisor of $\bar{H}+\alpha v$. In particular this implies that there is a (with positive radius) $d$-dim disc $\mathbb{B}\subseteq [0,1]^q$ with $\bar{H}$ as its centre such that $K_{U_t}+H$ is $\epsilon-\beta$-lc and  $R_H\subseteq R_{H_\alpha}$ for any $H\in \mathbb{B}$. This is a contradiction with the way we chose the sequence $\{H_1,\dots, H_m,\dots\}$. The function $R:\mathcal{H}_t \longrightarrow \mathbb{N}$ gives a finite decomposition of the set $\mathcal{H}_t$. This means that there are only a finite number of partial resolutions for all $(U_t,H)$ where $H\in \mathcal{H}_t$ for a fixed $t$. Using Noetherian induction completes the proof. 
$\Box$

\end{proof}

Now we prove a statement similar to [Sh2, 4.2].

\begin{lem}
Let $\mathcal{U}=\{(U, \Supp D)\}$ be a bounded family where we assume that each $(U, D)$ is Klt and exceptional and $K+D$ is antinef. Then the singularity is bounded i.e. there is a constant $\gamma >0$ such that each $(U,D)$ is $\gamma$-lc.   
\end{lem}
\begin{proof}

For $(U, \Supp B_i)$ a member of the family let 

\[\mathcal{H}_i = \{H =\sum h_{k,i}D_{k,i} |~ K + H ~is~ log~ canonical~ and~   -(K + H)~ is~ nef\}\]

where $D_i=\sum d_{k,i}D_{k,i}$.

It is a closed subset of a multi-dimensional cube (with bounded dimension) and so it is compact. Let $a_i=\inf\{\mld(U_i,H)~:~H\in \mathcal{H}_i\}>0$. Since the family is bounded then $\{a_i\}$ is bounded from below.

$\Box$

\end{proof}

Now we return to the division in \ref{3-5-J} and deal with each case as follows:

\item (First case in \ref{3-5-J}) Perform a hat of the first type for $U:=V_1$ and $K_U+G_U:=K_{V_1}+D_1$ (so we blow up $E,E'$). Then we get $V_2:=U'$ and $K_{V_2}+D_2:=K_{U'}+G_{U'}$ as defined above and $Y_1:=Y_U$. Now $V_2$ would be as in step 11, 12 or 13 so we can perform the appropriate operation as explained in each case. If $V_2$ is as in step 11 then  $a(E,V_2,\Delta_{V_2})=1$ and $a(E',V_2,\Delta_{V_2})=1$. If $V_2$ is as in step 12 then $E$ or $E'$ is not exceptional so we have $a(E,V_2,\Delta_{V_2})=1$ or $a(E',V_2,\Delta_{V_2})=1$. But if $V_2$ is as in step 13 then we get $U_1$ as defined in step 13 and so $a(E,U_1,\Delta_{U_1})=1$ or $a(E',U_1,\Delta_{U_1})=1$. In the later case we define (replace) $(V_2,D_2):=(U_1,G_{U_1})$.

So whatever case we have for $V_2$ we have $a(A,V_2,\Delta_{V_2})=1$ at least for one $A\in exc(Y/X)$.

\item (Second case in \ref{3-5-J}) Here we perform a hat of second type for $U:=V_1$ and $K_U+G_U:=K_{V_1}+D_1$ to get $V_2:=U'$ and $K_{V_2}+D_2:=K_{U'}+G_{U'}$. If $V_2$ is as in step 11 then  $a(E,V_2,\Delta_{V_2})=1$. If $V_2$ is as in step 12 then go to step 17. But if $V_2$ is as in step 13 then we get $U_1$ as defined in step 13 where $K_U+G_U:=K_{V_2}+D_2$ and then continue the process for $U_1$ as in step 15. 

Here in some cases we may not be able to make the singularities better for $K+\Delta$ immediately on $V_2$ but the algorithm ensures us that we will be able to do that in later steps.

\item (Third case in \ref{3-5-J}) In this case $V_1$ is 2-bir. We perform a hat of the third type where $U:=V_1$ and $K_U+G_U:=K_{V_1}+D_1$ so we get $V_2:=U'$ and $Y_1:=Y_U$ and $K_{V_2}+D_2:=K_{U'}+G_{U'}$. If $V_2$ is as in step 11 then  $a(E,V_2,\Delta_{V_2})=1$ and $a(C,V_2,\Delta_{V_2})=1$ ($E$ is the blown divisor and $C$ is on $V_1$, as in step 12 for $U:=V_1$). If $V_2$ is as in step 12 then $a(E,V_2,\Delta_{V_2})=1$ or $a(C,V_2,\Delta_{V_2})=1$. Now if $V_2$ is as in step 13 then we get $U_1$ as defined in step 13 and so $a(E,U_1,\Delta_{U_1})=1$ or $a(C,U_1,\Delta_{U_1})=1$. Then we define (replace) $(V_2,D_2):=(U_1,G_{U_1})$. 

 So whatever case we have for $V_2$ we have $a(A,V_2,\Delta_{2})=1$, after the appropriate operations, at least for one $A\in exc(Y/X)$.

\item After finitely many steps we get $V_r$ where $W/V_r$ such that $K_W+D:={^*(K_{V_r}+D_r)}$ with effective $D$ where $V_r$ is bounded. Since all the coefficients of $B_{V_r}$ are $\geq \frac{1}{2}-\tau$ ($B_{V_r}$ is the birational transform of $B_W$ where $K_W+B_W={^*K_X}$) then $(V_r, B_{V_r})$ is also bounded by lemma \ref{bound-1}. By construction $\Supp D_r=\Supp B_{V_r}$ and so $(V_r, D_r)$ is bounded. Lemma \ref{mainlemma} implies the boundedness of $W$ and so of $X$.

$$  \xymatrix{
    W \ar[d]\ar[r] & W \ar[d]\ar[r] & W \ar[d]\ar[r] & \dots\ar[r] & W\ar[d]&\\
    Y \ar[d]^{f}\ar[rd]^{g}  &  Y_1\ar[d]^{f_1}\ar[rd]^{g_1} & Y_2 \ar[d]\ar[rd] &\dots & Y_{r-1}\ar[d]\ar[rd]&\\
    X & V_1 & V_2 & \dots & V_{r-1}& V_r\\
}$$

\end{enumerate}
$\Box$
\end{proof}

\begin{cor}\label{BAB'}
The BAB Conjecture(\ref{BAB}) holds in dim 2.
\end{cor}
\begin{proof}
  {\textbf{ Reduction to the case $B=0$:}} We run the anti-LMMP on the divisor $K_X$; if there is an extremal ray $R$ such that $K_X.R>0$ then contract $R$ to get $X\longrightarrow X_1$. Note that $B.R<0$ so the bigness of $-K_X$ will be preserved (So $R$ has to be of birational type). Repeat the same process for $X_1$ i.e. if there is an extremal ray $R_1$ such that $K_{X_1}.R_1>0$ then contract it and so on. Since in each step we get a Pseudo-WLF then the canonical class can not become nef. Let $\overline{X}$ be the last model in our process, then $-K_{\overline{X}}$ is nef ad big. Now the boundedness of $\overline{X}$ implies the boundedness of $X$. So we replace $X$ by $\overline{X}$ i.e. from now on we can assume $B=0$. 

By theorem \ref{weak-2dim} $(X,0)$ has an $(\epsilon,n)$-complement $K_X+B^+$ for some $n\in \mathcal{N}_{\delta},2,\{0\}$. Now let $W\longrightarrow X$ be a minimal resolution and $\phi:W\longrightarrow S$ be the map obtained by running the classical MMP on $W$ i.e. contracting $-1$-curves to get a minimal $S$. As it is well known $S$ is $\mathbb{P}^2$ or a smooth ruled surface with no $-1$-curves. 

Let  $B_{S}^+=\sum b_{i,S}^+B_{i,S}^+$ be the pushdown of $B_{W}^+$ on $S$ where $K_W+B_{W}^+$ is the crepant pullback of $K_X+B^+$. Then define 

\[A_S:=\frac{b_{1,S}^+}{2}B_{1,S}^++ \sum_{i\neq 1} b_{i,S}^+B_{i,S}^+\]

If $S=\mathbb{P}^2$ then $-(K_S+A_S)$ is ample and $\Supp A_S=\Supp B_{S}^+$. By lemma \ref{bound-1} $(S, \Supp A_S=\Supp B_{S}^+)$ is  bounded. Then lemma \ref{mainlemma} implies the boundedness of $W$ and so of $X$.  

Now assume that $S$ is a ruled surface. If there is no exceptional curve (with negative self-intersection) on $S$ then $-(K_S+A_S)$ is nef and big if we take $B_{1,S}^+$ a non-fibre component of $B_{S}^+$. Since $S$ is smooth then $S$ is bounded and so $(S, \Supp A_S=\Supp B_{S}^+)$ is  bounded. 

But if there is an exceptional divisor $C$ on $S$ then contract $C$ as $S\longrightarrow S'$. So $S$ is a minimal resolution of $S'$. Since $\rho(S)=2$ and $(S',0)$ is $\delta$-lc then the index of each integral divisor on $S'$ is bounded. So $S'$ is bounded and then $(S',\Supp B_{S'}^+)$ is also bounded. This implies the boundedness of $S$, $W$ and so of $X$. Note that $B_{S'}^+\neq 0$ as $S'$ is WLF. 

$\Box$
\end{proof}

\begin{cor}\label{weak-proof}
Conjecture \ref{weak} holds for any finite set $\Gamma_f\subseteq [0,1]$ of rational numbers.
\end{cor}
\begin{proof}
It follows from corollary \ref{BAB'} 
$\Box$.
\end{proof}

%%%%%%%%%%%%%%%%%%%%%%%%%%%%%%%%%%%%%%%%%%%%%%%%%%%%%%%%%%%%%%%%%%%

\subsection{Second proof of the global case}

Remember that all the varieties are algebraic surfaces unless otherwise stated. We first prove the boundedness of varieties and then prove the boundedness of complements. This is somehow the opposite of what we did in the last subsection. However our proof was inspired by the theory of complements. The following proof heavily uses the properties of surfaces. That means that it is not expected to have a higher dimensional generalisation. The method also has some similarity with the proof of Alexeev and Mori [AM] in the sense that we both analyse a series of blow ups, but in different ways.

\begin{thm}\label{3-5-H}
The $BAB_{\delta, 2,[0,1]}$ holds.
\end{thm}

\begin{proof}
Now we reduce to the case $B=0$. Run the anti LMMP on the pair $(X, 0)$ i.e. if $-K_{X}$ is not nef then contract an extremal ray $R$ where $K_{X}.R>0$. This obviously contracts a curve in $B$. Repeating this process gives us a model $(X', 0)$ where $-K_{X'}$ is nef and big. Otherwise $X'$ should be with Picard number one and $K_{X'}$ nef. But this is impossible by our assumptions. We will prove the boundedness of $\{X'\}$ and so it implies the boundedness of $\{X\}$. Now we replace $(X,B)$ with $(X',0)$ but we denote it as $(X,0)$. We also assume that $\delta<1$ otherwise $X$ will be smooth and so with bounded index. 

Let $W\longrightarrow X$  be a minimal resolution. The main idea is to
prove that  there are only a bounded number of possibilities for the coefficients in $B_{W}$ where $K_{W}+B_{W}={^*K_{X}}$ i.e. the index of $K_{X}$ is bounded. 

{\bf{Strategy:}} Here we again have the familiar division into non-exceptional and exceptional cases.

First assume that $(X,0)$ is non-exceptional. So there will be a  $(0,n)$-complement $K_{X}+B^+$ for $n<58$. If we run the classical MMP on the pair $(W, 0)$ then we end up with $S$ which is either $\mathbb{P}^2$ or a ruled surface. Since $-(K_{S}+B_{S})={-_*(K_{W}+B_{W})}$ is nef and big then $K_{S}$ can not be nef. Let $K_{W}+B_{W}^+={^*(K_{X}+B_{X}^+)}$

\begin{lem}\label{3-5-D}
Let $G$ be a component in the boundary $B_{S}^+$ where $K_{S}+B_{S}^+={_*(K_{W}+B_{W}^+)}$ then $G^2$ is bounded from below and above. Moreover
there are only a bounded number of components in $B_{S}^+$.
\end{lem}
\begin{proof}
The boundedness of $G^2$ follows from the next lemma and the fact that $X$ is $\delta$-lc. The boundedness of number of components in $B_{S}^+$ is left to the reader. 
$\Box$
\end{proof}

The more general lemma below will also be needed later. 

\begin{lem}\label{3-5-E}
Let $(T,B_{T})$ be an $\delta$-lc WLF pair where $T$ is either $\mathbb{P}^2$ or a smooth ruled surface (with no $-1$-curves) and suppose $K_{T}+\overline{B}$ is antinef (lc) for a boundary $\overline{B}$. Let $M, B'_{T}$ be effective divisors with no common component such that $\overline{B}=B'_{T}+M$ then $M^2$ is bounded from above. 
\end{lem}

\begin{proof}
First assume that $T=\mathbb{P}^2$. In this case the lemma is obvious because if $M^2$ is too big then so is $\deg M$ and so it contradicts the fact that $\deg M \leq 3$. 

Now assume that $T$ is a ruled surface where $F$ is a general fibre other than those curves in the boundary and $C$ a section. The Mori cone of $T$ is generated by its two edges. $F$ generates one of the edges. If all the components of $M$ are fibres then $M^2=0$ and we are done. So assume otherwise and let $M\equiv aC+bF$ then $0<M.F=(aC+bF).F=a$ so $a$ is positive. Let $C^2=-e$ and consider the following two
cases:

1. $e\geq 0$:  We know that $K_{T}\equiv -2C+(2g-2-e)F$ where $g$ is a non-negative number [H, V, 2.11]. So we have 

\[ 0\geq (K_{T}+M+tC).F= -2+a+t \]
for some $t\geq 0$ where $B'_{T}\equiv tC+uF$ ($u\geq 0$ since $e\geq 0$). Hence $a+t\leq 2$. Calculations give $M^2=a(2b-ae)$. Since $a$ and $e$ are both non-negative then $M^2$ big implies that $b$ is big. But on other hand we
have:

\[ 0\geq (K_{T}+M+tC).C= (-2+a+t)(-e) +2g-2-e+b\]
This gives a contradiction if $b$ is too big because $e$ is also bounded. The boundedness of $e$ follows from the fact that $T$ is $\delta$-lc. In order in the local isomorphic subsection we proved that exceptional divisors have bounded selfintersection numbers. 

2. $e<0$: in this case by [H, V, 2.12] we have $e+2g\geq 0$ and so: 

\[ 0\geq (K_{T}+M).C=\]
 \[(-2+a)(-e) +2g-2-e+b=2g+e-2-(ae/2)+ (2b-ae)/2\]
Now since $2g+e-(ae/2)\geq 0$ then $(2b-ae)/2\leq 2$. So $M^2$ is bounded because $a$ is also bounded.
$\Box$
\end{proof}

 Let $P\in X$ be a singular point. If $P$ is not in the support of $B^+$ then the index of $K_{X}$ at $P$ is at most $57$ and so bounded. Now suppose that $P$ is in the support of $B^+$. If the singularity of $P$ is of type $E_{6}$, $E_{7}$, $E_{8}$ or $D_{r}$ then again the index of $K_{X}$ at $P$ is bounded. So assume that the singularity at $P$ is of type $A_{r}$. The goal is to prove that the number of curves in $exc(W/P)$ is bounded. In order we should prove that the number of $-2$-curves are bounded because the number of other curves is bounded by the proof of local isomorphic case. Note that the coefficient of any $E\in exc(W/P)$ in $B_{W}^+$ is positive and there are only a bounded number of possibilities for these coefficients. Let $\mathcal{C}$ be the longest connected subchain of $-2$-curves in $exc(W/P)$. Run the classical MMP on $W$ to get a model $W'$ such that there is a $-1$-curve $F$ on $W'$ s.t. it is the first $-1$-curve that intersects the chain $\mathcal{C}$ (if there is no such $W'$ and $F$ then $\mathcal{C}$ should consist of a single curve). We have two cases:

1. $F$ intersects, transversally and in one point, only one curve in $\mathcal{C}$ say $E$. First suppose that $E$ is a middle curve i.e. there are $E'$ and $E''$ in the
chain which both intersect $E$. Now contract $F$ so $E$ becomes a $-1$-curve. Then contract $E$ and then $E'$ and then all those which are on the side of $E'$. In
this case by contracting each curve we increase $E''^2$ by one. And so $E''$ will be a divisor on $S$ in $B_{S}^+$ with high self-intersection. By the lemma above there can
be only a bounded number of curves in $\mathcal{C}$ on the side of $E'$. Similarly there are only a bounded number of curves on the side of $E''$. So we are done in
this case. 

Now suppose that $E$ is on the edge of the chain and intersects $E'$. Let $t_{E}$ and $t_{F}$ be the coefficients of $E$ and $F$ in $B_{W}^+$ and similarly for
other curves. Let $h$ be the intersection number of $F$ with the curves in $B_{W'}^+$ except those in $\mathcal{C}$ and $F$ itself. Now we have 
\[ 0=(K_{W'}+B_{W'}^+).F= t_{E}+h-1-t_{F}\]
So  $h=1+t_{F}-t_{E}$. If $h\neq 0$ then $F$ intersects some other curve not in the chain $\mathcal{C}$. By contracting $F$ then $E$ and then other curves in the chain we get a
contradiction again. Now suppose $h=0$ i.e. $t_{E}=1$ and $t_{F}=0$. In this case let $x$ be the intersection of $E$ with the curves in $B_{W'}^+$ except those
in $\mathcal{C}$. So we have 
\[ 0=(K_{W'}+B_{W'}^+).E= -2t_{E}+t_{E'}+x\]
So $x=2t_{E}-t_{E'}>0$ and similarly we again get a contradiction.

2. Now assume that $F$ intersects the chain in more than one curve or intersects a curve with intersection number more than one. Suppose the chain $\mathcal{C}$ consists of
$E_{1}, \dots, E_{s}$ and $F$ intersects $E_{j_{1}}, \dots, E_{j_{l}}$.  Note that $l$ is bounded. If $F.E_{j_{k}}>1$ for all $0\leq k \leq l$ then contract $F$. So $E_{j_{k}}^2\geq 0$ after contraction of $F$ and they will not be contracted later and so they appear in the boundary $B_{S}^+$.
Now replace $\mathcal{C}$ with longest connected subchain when we disregard all $E_{j_{k}}$. Now go to step one again and if it doesn't hold come back to step two and so on. 

Now suppose $F.E_{j_{k}}=1$ for some $k$. So $F$ should intersect at least another $E_{j_{t}}$ where $t=k+1$ or $t=k-1$. Now contract $F$ so
$E_{j_{k}}$ becomes a $-1$-curve and would intersect $E_{j_{t}}$. Contracting  $E_{j_{k}}$ and possible subsequent $-1$-curves will prove that there are  a bounded number of curves between  $E_{j_{t}}$ and $E_{j_{k}}$. Now after contracting $E_{j_{k}}$ and all other curves between $E_{j_{t}}$ and $E_{j_{k}}$
 we will have $E_{j_{m}}^2\geq 0$ for each $m\neq k$. So again we take the longest connected subchain excluding all  $E_{j_{t}}$. And repeat the procedure. It
should stop after a bounded number of steps because the number of curves in $B_{S}^+$ is bounded. 
This boundedness implies that there are only a bounded number of possibilities for the coefficients in $B_{W}$ where $K_{W}+B_{W}={^*K_{X}}$. By
Borisov-McKernan $W$ belongs to a bounded family and so complements would be bounded. 

Here the proof of the non-exceptional case finishes and from now on we assume that $(X,0)$ is exceptional.

 Let $W\longrightarrow X$  be a minimal resolution. Let $0<\tau<\frac{1}{2}$ be a  number and the minimal log discrepancy of $(X,0)$ be $a=\mld(X,0)$. If $a\geq \frac{1}{2}+\tau$ then we know that $X$ belongs to a bounded family  according to step 1 in the proof of theorem \ref{weak-2dim} above. So we assume $a<\frac{1}{2}+\tau$ and then blow up all exceptional/$X$ curves $E$ with log discrepancy $a_{E}=a(E,X,0)\leq \frac{1}{2}+\tau$ to get $Y\longrightarrow X$ and put $K_{Y}+B_{Y}={^*K_{X}}$. Fix $E_1$, one of these exceptional divisors. Let $t\geq 0$ be a number such that there is an extremal ray $R$ such that $(K_{Y}+B_{Y}+tE_1).R=0$ and $E_1.R>0$ (and s.t.  $K_{Y}+B_{Y}+tE_1$ is Klt and antinef). Such $R$ exists otherwise there is a $t>0$ such that $K_{Y}+B_{Y}+tE_1$ is lc (and not Klt) and antinef. This is a contradiction  by [Sh2, 2.3.1]. Now contract $R: Y\longrightarrow Y_1$ if it is of birational type.

Again by increasing $t$ there will be an extremal ray $R_1$ on $Y_1$ such that $(K_{Y_1}+B_{Y_1}+tE_1).R_1=0$ and $E_1.R_1>0$ (preseving the nefness of $-(K_{Y_1}+B_{Y_1}+tE_1)$ ). If it is of birational type then contract it and so on. After finitely many steps we get a model $(V_1, B_{V_1}+t_1E_1)$ and a number $t_1>0$ with the following possible outcomes:

\begin{equation}\label{division}
\end{equation}
\begin{description}
 \item[$\diamond$] $ (V_1, B_{V_1}+t_1 E_1)$ is Klt, $\rho(V_1)=1$ and $K_{V_1}+B_{V_1}+t_1 E_1\equiv 0$.
 \item[$\diamond$] $ (V_1, B_{V_1}+t_1 E_1)$ is Klt and $\rho(V_1)=2$ and there is a non-birational extremal ray $R$ on $V_1$ such that $(K_{V_1}+B_{V_1}+t_1 E_1).R=0$. Moreover $K_{V_1}+B_{V_1}+t_1 E_1$ is antinef.
\end{description}

Note that for each element $E\in exc(Y/X)$, either $E$ is a divisor on $V_1$ or it is contracted to a point in the support of $E_1$.

\begin{lem}\label{3-5-C} 
For any $h>0$ there is an $\eta>0$ such that if $(T, B)$ is a $\delta$-lc pair ($\delta$ is already fixed) with a component $C$ of $B$ passing through $P\in T$, with a coefficient $t \geq h$, then either $K_{T}$ is $\delta +\eta$-lc at $P$ or $1-a_{E}>\eta$ for each exceptional divisor $E/P$ on a minimal resolution of $T$ ($a_E=$ log discrepancy of $(T,B)$ at $E$).
\end{lem}

\begin{proof}
 If $P$ is smooth or has $E_{6}, E_{7}$, $E_{8}$ or $D_{r}$ type of singularity then the lemma is clear since the index of $K_T$ at $P$ is bounded in all these cases (see the local isomorphism subsection). In order in all these cases there will be an $\eta>0$ such that $K_T$ is $\delta +\eta$-lc at $P$.

Now suppose that the singularity at $P$ is of type $A_{r}$. Take a minimal resolution  $W_T\longrightarrow T$ with $exc(W_T/P)=\{E_1, \dots, E_r\}$ (notation as in the local isomorphic subsection) and suppose that $j$ is the maximal number such that $\mld(P,T,0)=a'_j$ ( here we show the log discrepancy of $(T,0)$ at $E_{*}$ as $a'_{*}$) for an exceptional divisor $E_{j}/P$. Actually we may assume that $r-j$ is bounded. By the local isomorphic subsection , the distance of $E_{j}$ from one of the edges of $exc(W_T/P)$ is bounded. We denote the birational transform of $C$ on $W_T$ again by $C$. Suppose $C$ intersects $E_k$ in $exc(W_T/P)$. If $k\neq 1$ or $r$ then we have $(-E_{k}^2)a_{k}-a_{k-1}-a_{k+1}+x=0$ where $a_{*}$ shows the log discrepancy of the pair $(T,B)$ at $E_{*}$ and $x\geq h$ a number. So either $a_{k-1}-a_{k}\geq \frac{h}{2}$ or $a_{k+1}-a_{k}\geq \frac{h}{2}$. In either case the distance of $E_{k}$ is bounded from one of the edges of
$exc(W_T/P)$. If this edge is the same edge as for $E{j}$ then again the lemma is clear since the coefficients of $E_{k}$ and $E_{j}$ in $^*C$ (now $C$ is on
$T$ and $^*C$ on $W_T$) are bounded from below (in other words they are not too small). Now assume the otherwise i.e. $E_{k}$ and $E_{j}$ are close to different edges.  In this case
we claim that the coefficients of the members of $exc(W_T/P)$ in $\overline{B}_{W_T}$, where $K_{W}+\overline{B}_{W_T}={^*(K_{T}+tC)}$, are bounded from below. Suppose that the smallest coefficient
occurs at $E_{m}$. Simple calculation shows that we can assume that $E_{m}$ is one of the edges of $exc(W_T/P)$. Hence $E_{m}$ is with a bounded distance from $E_{j}$ or from $E_{k}$. 

Suppose that $E_{m}$ is with a bounded distance from $E_{j}$. If $a'_j\geq \frac{1+\delta}{2}$ then $K_T$ is $\frac{1+\delta}{2}$-lc at $P$. So we can assume that $a'_j<\frac{1+\delta}{2}$. We prove that all the numbers $1-a'_j, \dots, 1-a'_r$ are bounded from  below. In order, if $1<j<r$ then $(-E_{j}^2)a'_{j}-a'_{j-1}-a'_{j+1}=0$ (note that $-E_{j}^2>2$ in this case). Now  if $a'_{j-1}-a'_{j}\geq \frac{\delta}{2}$ then the chain will be bounded and so the index of $K_T$ at $P$. But if $a'_{j+1}-a'_{j}\geq \frac{\delta}{2}$ then $a'_{r}-a'_{r-1}\geq \frac{\delta}{2}$ and so $(-E_{r}^2-1)a'_{r}=1-(a'_{r}-a'_{r-1})\leq 1-\frac{\delta}{2}$. Hence if $m=r$ then we are done. But if $m=1$ then again the whole chain is bounded and so the index of $K_T$ at $P$. Now if $j=r$ then again the chain is bounded if $m=1$ and $a'_m=a'_j=a'_r<\frac{1+\delta}{2}$ if $m=r$.

 In the second case i.e. if $E_m$ is with a bounded distance from $E_{k}$ then the coefficient of $E_m$ in $^*C$ on $W$ is bounded from below. 
$\Box$
\end{proof}

\begin{lem}\label{3-5-I} 
For any $h>0$ there is a $\gamma>0$ such that if $(T, B)$ is a $\delta$-lc pair ($\delta$ is already fixed), WLF  with a component $C$ of $B$ passing through $P\in T$ and $t\geq h$ where $t$ is the coefficient of $C$ in $B$, then  $K_{T}$ is $\delta +\gamma$-lc. 
\end{lem}
\begin{proof}
 As discussed in lemma \ref{3-5-C} we may assume that the singularity at $P$ is of type $A_r$ and $1-a_{k}>\eta$ for some fixed number $\eta>0$ where $a_k$ is the log discrepancy of the pair $(T,B)$ at any exceptional divisor $E_{k}/P$ on $W_T$ where $W_T\longrightarrow T$ is a minimal resolution (we put $exc(W_T/P)=\{E_1, \dots, E_r\}$). Let $\mathcal{C}$ be the longest connected sub-chain of
$-2$-curves in $exc(W_T/P)$ and $W_{1}$  a model where $\mathcal{C}$ is intersected by a $-1$-curve $F$ for the first time i.e. we blow down $-1$-curves on $W_T$ till we get a model $W_1$ and a morphism $W_T\longrightarrow W_{1}$ such that $W_1$ is the first model where there is a $-1$-curve $F$ intersecting $\mathcal{C}$ (on $W_1$). Let $K_{W_T}+{B^+}\equiv 0$ be a (lc) $\mathbb{Q}$-complement of $K_{W_T}+B_{W_T}$. Assume that $F$ intersects $E_{j}$ in $\mathcal{C}$ and let $t_{E_{j}}$  and $t_{F}$ be the coefficients of $E_{j}$ and $F$ in $B^+$ on $W_T$ (similar notation for the coefficients of other exceptional divisors). Then an argument as in the proof of the non-exceptional case gives a contradiction: 

1. Suppose $F$ intersects, transversally and in one point, only one curve in $\mathcal{C}$ (which is $E_j$). First suppose that $E_{j}$ is a middle curve i.e. there are $E_{j-1}$ and $E_{j+1}$ in  $\mathcal{C}$  which both intersect $E_{j}$. Now contract $F$ so $E_{j}$ becomes a $-1$-curve. Then contract $E_{j}$ and then $E_{j-1}$ and then all those which are on the of $E_{j-1}$. By contracting each curve we increase $E_{j+1}^2$ by one. If we continue contracting $-1$-curves we get $S$ ($S=\mathbb{P}^2$ or a ruled surface with no $-1$-curve) where $E_{j+1}$ is a component of $B_{S}$. By lemma \ref{3-5-E} there can be only a bounded number of curves in $\mathcal{C}$ on the side of $E_{j-1}$. Similarly there are only a bounded number of curves in $\mathcal{C}$ on the side of $E_{j+1}$. So we are done in this case. 

Now suppose that $E_{j}$ is on the edge of the chain $\mathcal{C}$ and it intersects $E_{j-1}$.  Let ${B^+}_{W_{1}}=\dot{B^+}+M$ ($M$ and $\dot{B^+}$ with no common component) where each component of $\dot{B^+}$ is either $F$ or an element of $\mathcal{C}$. Now we have 
\[ 0= (K_{W_{1}}+{B^+}_{W_{1}}).F= t_{E_{j}}-1-t_{F}+(M.F)\]
So $M.F= 1+t_{F}-t_{E_{j}}$. Similarly let ${B^+}_{W_{1}}=\ddot{B^+}+N$ ( $N$ and $\ddot{B^+}$ with no common component) where each component of $\ddot{B^+}$ is either $F$ or an element of $\mathcal{C}$. Then we have
\[ 0= (K_{W_{1}}+{B^+}_{W_{1}}).E_{j}= -2t_{E_{j}}+t_{E_{j-1}}+t_{F}+(N.E_{j})\]
 and so $t_{E_{j}}=t_{E_{j-1}}-t_{E_{j}}+t_{F}+(N.E_{j})>\eta$. Hence $ t_{E_{j-1}}-t_{E_{j}}>\frac{\eta}{3}$ or $t_{F}>\frac{\eta}{3}$ or $(N.E_{j})>\frac{\eta}{3}$. 

If $t_{F}>\frac{\eta}{3}$ then by contracting $F$ we increase $M^2$ at least by $(M.F)^2\geq t_{F}^2> (\frac{\eta}{3})^2$. We have the same increase when we contract $E_{j}$ and then $E_{j-1}$ and so on. So  lemma \ref{3-5-E} shows the boundedness of $\mathcal{C}$. 

If $(N.E_{j})>\frac{\eta}{3}$ then proceed similar to the last paragraph.

 If $ t_{E_{j-1}}-t_{E_{j}}>\frac{\eta}{3}$ then $ t_{E_{j-1}}>t_{E_{j}}+\frac{\eta}{3}$. This implies that $t_{E_{j}}\leq 1-\frac{\eta}{3}$ then $M.F\geq \frac{\eta}{3}$ and so we continue as above. 

2. Now assume that $F$ intersects $\mathcal{C}$ in more than one curve or intersects a curve in   $\mathcal{C}$ with intersection number more than one. Suppose the chain $\mathcal{C}$ consists of
$E_{s}, \dots, E_{u}$ and $F$ intersects $E_{j_{1}}, \dots, E_{j_{l}}$. Note that $l$ is bounded.

If $F.E_{j_{k}}>1$ for all $1\leq k \leq l$ then contract $F$. So  $E_{j_{k}}^2\geq 0$ after contraction of $F$ and hence $E_{j_{k}}$ can  not be contracted  and so it appears in the boundary on a ``minimal'' model $S$ (i.e. $S$ is the projective plane or a smooth ruled surface with no $-1$-curve).
 Replace $\mathcal{C}$ with its longest connected subchain when we disregard all $E_{j_{k}}$. From here we can go back to step one and repeat the argument. 

Now suppose $F.E_{j_{k}}=1$ for some $k$. So $F$ should intersect at least another $E_{j_{q}}$ where $q=k+1$ or $q=k-1$. Now contract $F$ so
$E_{j_{k}}$ becomes a $-1$-curve and would intersect $E_{j_{q}}$. Contracting $E_{j_{k}}$ and possible subsequent $-1$-curves will prove that there are only a bounded number of curves between $E_{j_{q}}$ and $E_{j_{k}}$ in $\mathcal{C}$. Now after contracting $E_{j_{k}}$ and all other curves between $E_{j_{q}}$ and $E_{j_{k}}$
 we will have $E_{j_{m}}^2\geq 0$ for each $m\neq k$. So again we take the longest connected subchain excluding $E_{j_{1}}, \dots,E_{j_{l}} $ and go back to step one.

 This process should stop after a bounded number of steps because the number of curves in $B_{S}^+$ with coefficient $>\eta$ is bounded ($S$ is again a ``minimal'' model). To prove this later boundedness note that $(K_{S}+{B^+}_{S}).F=0$, where we assume that $S$ is a ruled surface and $F$ a fibre, implies that there are only a bounded number of non-fibre components in $B_{S}^+$ with coefficient $>\eta$. Let $L$ be a section and $t_L$ be its coefficient in $B_{S}^+$ and $F_i$ fibre components of  $B_{S}^+$ with $t_{F_i}>\eta$. So 

\[0 \geq (K_{S}+t_{L}L+\sum_{i}t_{F_i}F_i).L= (-2L+(2g-2-e)F+t_{L}L+\sum_{i}t_{F_i}F_i).L\]\[=-t_Le+e+2g-2+\sum_{i}t_{F_i}\] 

which proves that there are a bounded number of $F_i$ ($L^2=-e$ and $e+2g\geq 0$ if $e<0$). So the chain $\mathcal{C}$ should have a bounded length. This implies that if we throw $C$ away in the boundary $B$ then the mld at $P$ will increase at least by $\gamma >0$ a fixed number (which doesn't depend on $P$ nor $T$). This proves the lemma. 
$\Box$
\end{proof}

Lemma \ref{3-5-I} settles the first case in \ref{division}  by deleting the boundary $B_{V_1}$. Now assume the second case in the division above in \ref{division}. Let $F$ be a general fibre of the contraction defined by the extremal ray $R$. If the other extremal ray of $V_1$ defines a birational map $V_1\longrightarrow Z$ (otherwise delete the boundary and use \ref{3-5-I} ) then let $H$ be the exceptional divisor of this contraction. 

If $K_{V_1}$ is antinef then again use \ref{3-5-I}. If $K_{V_1}$ is not antinef  and if $E_1\neq H$ then apply lemma \ref{3-5-I} to $(Z,B_Z)$. Boundedness of $Z$ implies the boundedness of $V_1$ and so we can apply lemma \ref{bound-1}.  But if $K_{V_1}$ is not antinef and $E_1=H$ then perform a hat of the third type as defined in the proof of theorem \ref{weak-2dim} where $(U,G_U):=(V_1, B_{V_1}+t_1E_1)$ and $V_2:=U'$. We can use lemma  \ref{3-5-I} on $V_2$ or after contracting a curve on $V_2$, in order, to get the boundedness of $V_2$. Boundedness of $V_2$ implies the boundedness of $V_1$.  
$\Box$
\end{proof}

\begin{cor}
 Conjecture $WC_{\delta, 2,\Gamma_f}$ holds in the global case where $\Gamma_f$ is a finite subset of rational numbers in $[0,1)$.
\end{cor}

\begin{proof}
Obvious by theorem \ref{3-5-H}.
\end{proof}

%%%%%%%%%%%%%%%%%%%%%%%%%%%%%%%%%%%%%%%%%%%%%%%%%%%%%%%%%%%%%%%%%%%%

\subsection{An example}

\begin{exa}
Let $m$ be a positive natural number. For any $1>\eta >0$ and any $\tau > 0$ there is a model $(X,0)$ (may not be global) satisfying the followings:

\begin{enumerate}
 \item $X$ is $\frac{1}{m}$-lc. 
 \item Suppose $Y \longrightarrow X$ is a partial resolution such that $K_{Y}+B_{Y}={^*K_{X}}$ is $\frac{1}{m}+\eta$-lc and $b_{i}>\frac{m-1}{m}-\eta$. Put $D=\sum \frac{m-1}{m}B_{i}$. 
 \item $K_{Y}+D$ is not $\frac{1}{m}+\tau$-lc.
\end{enumerate}
\end{exa}
  
\begin{proof} Let $P\in X$ and $X$ smooth outside $P$. Suppose the minimal resolution of $P$ has the following diagram:

\begin{displaymath}
\xymatrix{O^{-3} \ar@{-}[r]& O^{-2} \ar@{-}[r]& \dots &\ar@{-}[r]& O^{-2}\ar@{-}[r]  &O^{-2}\ar@{-}[r] & O^{-4}}
\end{displaymath}
where the numbers show the self-intersections.

This diagram has the following corresponding system on a minimal resolution where $a_{i}$ stand for the log discrepancies:
\[ 3a_{1}-a_{2}-1=0\]
\[2a_{2}-a_{1}-a_{3}=0\]
\[\vdots\]
\[2a_{r-1}-a_{r-2}-a_{r}=0\]
\[4a_{r}-a_{r-1}-1=0\]

Now put $a_{r-1}-a_{r}=t$ so $a_{r-2}-a_{r-1}=t$, $\dots$, $a_{1}-a_{2}=t$. So we have: $a_{r}=\frac{1+t}{3}$ and $a_{1}=\frac{1-t}{2}$. The longer the chain is the smaller the $t$ is and the discrepancies vary from $-\frac{1+t}{2}$ to $\frac{t-2}{3}$. Other $a_{i}$ can be calculated as $a_{i}=a_{1}-(i-1)t=\frac{1-t}{2}-(i-1)t=\frac{1-(2i-1)t}{2}$.  

Suppose $j$ is such that $a_{j}<\frac{1}{m}+\eta$ but $a_{j-1}\geq \frac{1}{m}+\eta$. So the exceptional divisors corresponding to $a_{r}, a_{r-1}, \dots, a_{j}$ will appear on $Y$ but others not. Now we try to compute the log discrepancies of the pair $(Y, D)$. In order the minimal resolution for $P\in X$ is also the minimal resolution for $Y$. But here just $E_{1}, \dots, E_{j-1}$ are exceptional/$Y$. The system for the new log discrepancies (for $(Y, D)$) is as follows:
\[3a'_{1}-a'_{2}-1=0\] 
\[2a'_{2}-a'_{1}-a'_{3}=0\]
\[\vdots\]
\[2a'_{j-2}-a'_{j-3}-a'_{j-1}=0\]
\[2a'_{j-1}-a'_{j-2}-\frac{1}{m}=0\]

 Again put $a'_{j-2}-a'_{j-1}=s$ so similarly we have $a'_{j-1}=\frac{1}{m}+s$ and $a'_{1}=\frac{1-s}{2}$. If $j$ is big (i.e. if t is small enough) then $s$ would be small and so $a'_{j-1}=\frac{1}{m}+s <\frac{1}{m}+\tau$. Hence $(Y, D)$ is not $\frac{1}{m}+\tau$-lc. 

$\Box$
\end{proof}

%%%%%%%%%%%%%%%%%%%%%%%%%%%%%%%%%%%%%%%%%%%%%%%%%%%%%%%%%%%%%%%%%

\subsection{Local cases revisited}

Using the methods in the proof of the global case, we give a new proof of the local cases. Here again by $/Z$ we mean $/P\in Z$ for a fixed $P$. The following is the main theorem in this subsection.

\begin{thm}\label{4-A'} Conjecture  $WC_{\delta, 2,{\Phi_{sm}}}$ holds in the local case i.e. when we have $\dim Z\geq 1$.
\end{thm}

\begin{proof} Our proof is similar to the non-exceptional global case. Here the pair $(X/Z,B)$ is a relative WLF surface log pair (i.e. $-(K_X+B)$ is nef and big/$Z$), where $(X,B)$ is $\delta$-lc and $B\in \Phi_{sm}$. Fix $P\in Z$. Then there exists a regular $(0,n)$-complement/$P\in Z$, $K+B^+$ for some $n\in\{1,2,3,4,6\}$ by [Sh2]. 

\begin{enumerate}
\item Remember the first step in the proof of theorem \ref{weak-2dim}. 

\item Remember definition \ref{D-tau} and lemma \ref{transform}. Let $m$ be the smallest number such that $\frac{1}{m}\leq \delta$. Let $h=\min\{\frac{k-1}{k}-\frac{u}{r!}>0\}_{1\leq k\leq m} $ where $u,k$ are natural numbers and $r=\max\{m,6\}$. Now choose a $\tau$ for $m$ as in lemma \ref{transform} such that $\tau <h$. 

Blow up one exceptional divisor $E/P$ via $f:Y\longrightarrow X$ such that the log discrepancy satisfies $\frac{1}{k}\leq a(E,X,B)\leq \frac{1}{k}+\tau$ for some $k$ (if such $E$ doesn't exist then go to step 1). The crepant log divisor $K_Y+B_Y$ is $\frac{1}{m}$-lc and so by the choice of $\tau$, $K_Y+D_{\tau}$ is also $\frac{1}{m}$-lc ($D_{\tau}$ is constructed for $B_Y$). Let $K_Y+B_{Y}^+$ be the crepant blow up of $K_X+B^+$. Then again by the way we chose $\tau$ we have $D_{\tau}\leq B_{Y}^+$. Now run the anti-LMMP/$P\in Z$ over $K_Y+D_{\tau}$ i.e. contract any birational type extremal ray $R$/$P\in Z$ such that $(K_Y+D_{\tau}).R> 0$. At the end we get a model $X_1$ with one of the following properties:

\begin{description}
  \item[$\diamond$]  $(K_{X_1}+D_{\tau})\equiv 0/P\in Z$ and $K_{X_1}+D_{\tau}$ is $\frac{1}{m}$-lc.
  \item[$\diamond$] $-(K_{X_1}+D_{\tau})$ is nef and big/$P\in Z$ and $K_{X_1}+D_{\tau}$ is $\frac{1}{m}$-lc.
\end{description}

where $K_{X_1}+D_{\tau}$ is the birational transform of $K_{Y}+D_{\tau}$ and let $g:Y\longrightarrow X_1$ be the corresponding morphism.

The nefness of $-(K_{X_1}+D_{\tau})$ comes from the fact that  $D_{\tau}\leq B_{1}^+$. And $K_{X_1}+D_{\tau}$ is $\frac{1}{m}$-lc by applying lemma \ref{transform}.

\item Whichever case occurs above, to construct a complement, it is enough to bound the index of $K_{X_1}+D_{\tau}/P$.   
\item Let $C$ be a curve contracted by $g:Y\longrightarrow X_1$. If $C$ is not a component of $B_Y$ then the log discrepancy of $C$ with respect to $K_{X_1}+B_{X_1}$ is at least 1 where $K_{X_1}+B_{X_1}$ is the birational transform of $K_{Y}+B_{Y}$. Moreover $g(C)\in \Supp B_{X_1}\neq \emptyset$. So the log discrepancy of $C$ with respect to $K_{X_1}$ is more than 1. This means that $C$ is not a divisor on a minimal resolution $W_1\longrightarrow X_1$. Let $W\longrightarrow X$ be a minimal resolution. Then there is a morphism $W\longrightarrow W_1$. Hence $exc(W_1/X_1)\subseteq exc(W/X)\cup \Supp (B=B_X)$. Now if $C\in exc(W/X)\cup \Supp B$ is contracted by $g$ then $a(C,X_1,D_{\tau})<a(C,X,B)$.      

\item Let $(X_1,B_1):=(X_1,D_{\tau})$ and repeat the process. In other words again we blow up one exceptional divisor $E$ via $f_1:Y_1\longrightarrow X_1$ such that the log discrepancy satisfies $\frac{1}{k}\leq a(E,X_1,B_1)\leq \frac{1}{k}+\tau$ for some natural number $k>1$. The crepant log divisor $K_{Y_1}+B_{1,Y_1}$ is $\frac{1}{m}$-lc and so by lemma \ref{transform} $K_{Y_1}+D_{1,\tau}$ is $\frac{1}{m}$-lc. Note that the point which is blown up on $X_1$ can not be smooth since $\tau <h$ as defined above. So according to the last step the blown up divisor $E$ is a member of $exc(W/X)\cup \Supp B$. Now we again run the anti-LMMP on $K_{Y_1}+D_{1,\tau}$ and proceed as in step 2.

$$  
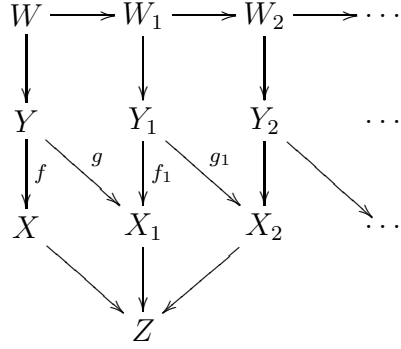
\xymatrix{
    W \ar[d]\ar[r] & W_1 \ar[d]\ar[r] & W_2 \ar[d]\ar[r] & \dots \\
    Y \ar[d]^{f}\ar[rd]^{g}  &  Y_1\ar[d]^{f_1}\ar[rd]^{g_1} & Y_2 \ar[d]\ar[rd] &\dots\\
    X \ar[rd]& X_1 \ar[d]& X_2 \ar[ld]& \dots \\
     &Z&& \\
}$$

\item Steps 4,5 show that each time we blow up a member of $exc(W/X)\cup \Supp B$ say $E$. And if we blow that divisor down in some step then the log discrepancy $a(E,X_j,B_j)$  will decrease. That divisor will not be blown up again unless the log discrepancy drops at least by $\frac{1}{2(m-1)}-\frac{1}{2m}$ (this is not a sharp bound). So after finitely many steps we get a model $X_i$ with a standard boundary $B_i$ such that  there is no $E/P$ where $\frac{1}{k}\leq a(E,X_i,B_i)\leq \frac{1}{k}+\tau$ for any $1<k\leq m$. Hence the index of $-(K_{X_i}+B_{i})/P$ is bounded and so we can construct an appropriate complement for $(X_i,B_i)/Z$. This implies the existence of the desired complement for $(X,B)/Z$.
\end{enumerate}
$\Box$
\end{proof}

\clearpage

%%%%%%%%%%%%%%%%%%%%%%
%%%%%%%%%%%%%%%%%%%%%%%%%%
%%%%%%%%%%%%%%%%%%%%%%%%

\section{$\epsilon$-log canonical  complements in higher dimensions}

 In this section we consider  the $(\epsilon, n)$-lc complements in higher dimensions i.e. in dimensions more than two. This is a joint work in progress with V.V. Shokurov.
 In subsection 2.1 we try to work out the proof of theorem \ref{weak-2dim} in dim 3 and we point out the problems we have to solve in order to finish the proof of conjecture \ref{weak} in dim 3 (this is the plan of the author). In subsection 2.2 we outline Shokurov's plan on the same problem. These plans have already won an EPSRC three years postdoctoral fellowship by the author. 

Let $X\longrightarrow Z$ be an extremal $K_X$-negative contraction where $X$ is a 2-dim Pseudo-WLF and $Z$ is a curve. We know that $Z\simeq \mathbb{P}^1$ since $Z$ should be rationally connected as $X$ is. Moreover $\rho(X)=2$. Similar Mori fibre spaces in higher dimensions are not that simple. This makes the boundedness problem of $(\epsilon, n)$-lc complements more difficult in higher dimensions. We also don't know yet whether the index of $K_X+B$ will be bounded if we fix the mld at a point. 

In section 1 we first proved the boundedness of $\epsilon$-lc complements and then the BAB. But in higher dimensions we expect to prove both problems together at once. In other words in some cases where it is difficult to prove the boundedness of varieties, it seems easier to prove the boundedness of complements;  specially when we deal with a fibre space. Conversely when it is difficult to prove the boundedness of $\epsilon$-lc complements, it is better to prove the boundedness of pairs; this is usually the case when the pairs are exceptional.

\begin{lem}\label{pre-2}
Let $X\dasharrow X'$ be a flip/$Z$ and assume that $(X,B)$ is $(\epsilon,n)$-complementary/$Z$ then $(X',B')$ is  $(\epsilon,n)$-complementary/$Z$ where $B'$ is the birational transform of $B$. 
\end{lem}
\begin{proof}
Obvious from the definition of $(\epsilon,n)$-complements.
\end{proof}

Note that in the previous lemma it doesn't matter that the flipping is with respect to which log divisor. 

\begin{lem}\label{pre-3}
Let $(Y,B)$ be a pair and $Y\dashrightarrow Y'/Z$ be a composition of divisorial contractions and flips$/Z$ such that in each step we contract an extremal ray $R$ where $(K+B).R\geq 0$. Suppose $B'=\sum {b'}_i{B'}_i$ is the birational transform of $B$, the pair $(Y',B')$ is $(\epsilon,n)$-complementary/$Z$ and $(n+1){b'}_i\geq n{b'}_i$ for each coefficient ${b'}_i$ then $(Y,B)$ is also $(\epsilon,n)$-complementary/$Z$.
\end{lem}
\begin{proof}
Clear by lemmas \ref{pre-1} and \ref{pre-2}.
\end{proof}

\begin{lem}  
The Klt Pseudo-WLF property is preserved under extremal flips and divisorial contractions with respect to any log divisor.
\end{lem}
\begin{proof} Let $X$ be a Klt Pseudo-WLF and $B$ a boundary such  that $(X,B)$ is a Klt WLF. Now let $X\dasharrow X'$ be an extremal flip corresponding to an extremal ray $R$. Since $(X,B)$ is a Klt WLF then there is a rational boundary $D$ such that $K_X+D$ is antiample and Klt. Now let $H'$ be an ample divisor on $X'$ and $H$ its transform on $X$. There is a rational $t>0$ such that $K_X+D+tH$ is antiample and Klt. Now take a Klt $\mathbb{Q}$-complement  $K_X+D+tH+A\equiv 0$. So we have $K_{X'}+D'+tH'+A'\equiv 0$ on $X'$. From the assumptions $K_{X'}+D'+A'$ is antiample and Klt. So $X'$ is also a Klt Pseudo-WLF. 

If $X\longrightarrow X'$ is a divisorial extremal contraction then proceed as in the flip case by taking an ample divisor $H'$ on $X'$.

$\Box$ 
\end{proof}

\begin{defn} Let $(V,B_V)$ and $(U,B_U)$ be  lc pairs. $U$ is called a semi-partial resolution of $V$ if there is a partial resolution $(W,B_W)$ of $(V,B_V)$ such that $W$ and $U$ are isomorphic in codim 1. 
\end{defn}

\begin{defn}[$D$-LMMP]\label{d-lmmp}
Let $D$ be an $\mathbb{R}$-Cartier divisor on a normal variety $X$. We say $D$-LMMP holds if the followings hold: 

\begin{description}
\item[$\diamond$] Any $D$-negative extremal ray $R$ on $X$ can be contracted. And the same holds in the subsequent steps for the birational transform of $D$.
\item[$\diamond$] If a $D$-contraction as in the first step is a flipping then the corresponding $D$-flip exists. 
\item[$\diamond$] Any sequence of $D$-flips terminates.
\end{description}
 
\end{defn}

If $D:=K+B$ for a lc $\mathbb{R}$-Cartier divisor $K+B$ then we know that $D$-LMMP holds in dim 3 by [Sh5]. 

\begin{rem}
Let $D$ be an $\mathbb{R}$-Cartier divisor on a variety $X$ of dim $d$ and assume that LMMP holds in dim $d$. Moreover assume that $\beta D\equiv K+B$ for a lc $\mathbb{R}$-Cartier log divisor $K+B$ and $\beta >0$ then the $D$-LMMP holds. Since in this case $D$-LMMP and $K+B$-LMMP are equivalent.
\end{rem}

\begin{exa}
Let $(X,B)$ be a $d$-dim Klt WLF and suppose LMMP holds in dim $d$ then $-K$-LMMP holds. In order since $(X,B)$ is a Klt WLF then there is a Klt $\mathbb{Q}$-complement $K+B^+\equiv 0$. There is a $t>0$ such that $K+B^++tB^+\equiv tB^+$ is Klt. Since $-K\equiv B^+$ then $-K$-LMMP is equivalent to $B^+$-LMMP and again equivalent to $tB^+$-LMMP. Since  $K+B^++tB^+$-LMMP holds so does $-K$-LMMP.

\end{exa}

%%%%%%%%%%%%%%%%%%%%%%%%%%%%%%%%%%%%%%%%%%%%%%%%%%%%%%%%%%%%%%%%%%%

\subsection{$\epsilon$-lc complements in dimension 3}

In this subsection we propose a plan toward the resolution of conjecture \ref{weak} in dim 3.

We repeat the proof of \ref{weak-2dim}, in dim 3, step by step:

\begin{enumerate}
\item Under the assumptions of conjecture \ref{weak} for $d=3$ and $\Gamma=\{0\}$, first assume that $(X,0)$ is non-exceptional.
\item We don't have much information about the accumulation points of mlds in dim 3. Actually we still have not proved ACC in dim 3. As pointed out in the introduction in section 1, only one case of ACC in dim 3 is remained to be proved. Remember that Shokurov's program tries to use complements in dim $d-1$ to prove the ACC in dim $d$. So it is reasonable to assume ACC in dim $d-1$. 

Lets show the set of accumulation points of mlds of $d$-dim lc pairs $(T,B)$, where $B\in \Gamma$, with $Accum_{d, \Gamma}$.

\item We need the inductive version of complements; since $(X,0)$ is not exceptional then it is expected that there is an inductive $(0,n)$-complement $K_X+B^+$ where $n \in \mathcal{N}_{2}$. 

\item Remeber definition \ref{D-tau}. We can similarly define $D_{\tau, A}$  for a boundary $B$, with respect to a real number $\tau\geq 0$ and a set $A\subseteq [0,1]$: 

\[D_{\tau, A}:=\sum_{b_{i}\notin [a-\tau, a]}b_{i}B_{i}+\sum_{b_{i}\in [a-\tau, a]} aB_{i}\] 

where in the first term  $b_{i}\notin [a-\tau, a]$ for any $a\in A$ but in the second term $a\in A$ is the biggest number satisfying $b_{i}\in [a-\tau, a]$. 

\begin{defn} Let $A\subseteq [0,1]$ and $(T,B)$ a log pair. We say that $(T,B)$ is $A$-lc if $(T,B)$ is $x$-lc where $x:=1-\sup\{A\}$. 
\end{defn}

Assuming the ACC in dim 3 a statement similar to lemma \ref{transform} may hold: For any $\gamma >0$ and finite set $A\subseteq [0,1]$ containing $1-\gamma$ there is a real number $\tau>0$ such that if $(T,B_T)$ is a 3-fold log pair, $P\in T$, $K_T+B_T$ is $\gamma$-lc in codim 2 at $P$ and $D_{\tau,A}\in A$  then  $K_T+D_{\tau,A}$ is also $\gamma$-lc in codim 2 at $P$. 

Moreover we expect to choose a $\tau>0$ such that the followings hold as well: 
\begin{itemize}
\item If $B_T\in A$ and $E$ the exceptional divisor of a smooth blow up of $T$ then $a(E,T,B_T)\notin [1-a, 1-a+\tau]$ for any $a\in A$.

\item  If $B_T\in A$ and the pair $(T,B_T)$ non-exceptional then we can refine $\mathcal{N}_2$ such that  there is a $(0,n)$-complement $K_T+B_{T}^+$ for some $n\in \mathcal{N}_2$ where $B_T\leq B_{T}^+$.
\end{itemize}

\item Let $A_1:=\{a_1\}$ where $1-a_1=\max Accum_{3,\{0\}}\cap [0,\delta]$. Now blow up all exceptional divisor $E$ such that $a(E,T,B_T)\in [1-a, 1-a+\tau]$ for some $a\in A_1$ to get $f:Y\longrightarrow X$. Construct $D_{\tau, A_1}$ for $B_Y$ where $K_Y+B_Y$ is the crepant pull back. So $(Y,D_{\tau, A_1})$ is $A_1$-lc. Run the $D$-LMMP where $D:=-(K_Y+D_{\tau, A_1})$. At the end we get $Y\dasharrow X_1$ and $X_1\dasharrow S_1$ such that $-(K_{X_1}+D_{\tau,A_1})$ is nef and $\equiv 0/S_1$ and $-(K_{S_1}+D_{\tau,A_1}).R>0$ for any birational type extremal ray $R$. 

\item There are the following possibilities for the model $S_1$:

\begin{description}

\item[$\diamond$] $\rho(S_1)=1$, $-(K_{S_1}+D_{\tau,A_1})=-(K_{S_1}+B^+)\equiv 0$ and $K_{S_1}+D_{\tau,A_1}$ is $A_1$-lc. 

\item[$\diamond$] There is a fibration type extremal ray $R$ such that , $-(K_{S_1}+D_{\tau,A_1}).R=0$ and $K_{S_1}+D_{\tau,A_1}$ is $A_1$-lc.

\item[$\diamond$]  $-(K_{S_1}+D_{\tau,A_1})$ is nef and big and $K_{S_1}+D_{\tau,A_1}$ is $A_1$-lc.
\end{description}

\item In the first case in the division above we are done. In the second and third case then replace $(X,0)$ by $(X_1,B_1):=(X_1,D_{\tau, A_1})$ and go back to step one and repeat. By repeating and repeating the process, each time we get new coefficients. In other words we need to replace $A_i$ with $A_{i+1}$ such that $A_i\subseteq A_{i+1}$. We need to prove that $\cup_{i\rightarrow\infty} A_i$ is finite.

\item At the end we get a model $(X_r, B_r)$ which is terminal in codim 2. Then we hope to prove the boundedness of the index of $K_{X_r}+B_r$ possibly after some more blow ups and blow downs. This will settle the problem if $-(K_{X_r}+B_r)$ is nef and big. Otherwise we may have a fibration and $K_{X_r}+B_{r}^+=K_{X_r}+B_r+N$ where $N$ is vertical. Then we may replace $N$ by $N'$ and construct a desirable complement  $K_{X_r}+B_r+N'$. At the end we need to prove that the boundedness of the complement implies the boundedness of the pairs. 

\item Now let $(X,0)$ be exceptional.  Since $BAB_{1,3,\{0\}}$ holds by [KMMT]
 then assuming ACC in dim 3, there is a $\tau>0$ such that $BAB_{1-\tau,3,\{0\}}$ also holds. Blow up an exceptional/$X$ divisor $E_1$ with log discrepancy $a_{E_1}=a(E_1,X,0)\leq 1-\tau$ to get $Y\longrightarrow X$ and put $K_{Y}+B_{Y}={^*K_{X}}$. Let $t\geq 0$ be a number such that there is an extremal ray $R$ such that $(K_{Y}+B_{Y}+tE_1).R=0$ and $E_1.R>0$ ( and s.t. $K_{Y}+B_{Y}+tE_1$ Klt and antinef). Such $R$ exists otherwise there is a $t>0$ such that $K_{Y}+B_{Y}+tE_1$ is lc (and not Klt) and antiample. This is a contradiction with the fact that $(X,0)$ is exceptional. Now contract $R: Y\longrightarrow Y_1$ if it is of birational type (and perform the flip if it is a flipping).

Again by increasing $t$ there will be an extremal ray $R_1$ on $Y_1$ such that $(K_{Y_1}+B_{Y_1}+tE_1).R_1=0$ and $E_1.R_1>0$ (preserving the nefness of $-(K_{Y_1}+B_{Y_1}+tE_1)$ ). If it is of birational type then contract it and so on. After finitely many steps we get a model $(V_1, B_{V_1}+t_1E_1)$ and a number $t_1>0$ with the following possible outcomes:

\begin{description}
 \item[$\diamond$] $ (V_1, B_{V_1}+t_1 E_1)$ is Klt, $\rho(V_1)=1$ and $K_{V_1}+B_{V_1}+t_1 E_1\equiv 0$.
 \item[$\diamond$] $ (V_1, B_{V_1}+t_1 E_1)$ is Klt  and there is a fibre type extremal ray $R$ on $V_1$ such that $(K_{V_1}+B_{V_1}+t_1 E_1).R=0$ and $K_{V_1}+B_{V_1}+t_1 E_1$ is antinef. 
\end{description}

If the second case occurs then we don't know $\rho(V_1)$ unlike the surface case where $\rho(V_1)=2$.

\item In the proof of theorem \ref{weak-2dim} we introduced three types of hat. Here also we can similarly define hats but it is not clear yet how to proceed.

\end{enumerate}

%%%%%%%%%%%%%%%%%%%%%%%%%%%%%%%%%%%%%%%%%%%%%%%%%%%%%%%%%%%%%%%%

\subsection{$\epsilon$-lc complements in dimension 3: Shokurov's approach}

Here we explain Shokurov's approach to the problem discussed in 4.1.

\begin{enumerate}

\item We know that the $BAB_{1,3,\{0\}}$ holds by [KMMT]. Let $a$ be the smallest positive real number with the following property: $BAB_{a',3,\{0\}}$ holds for any $a'>a$. The idea is to prove that $BAB_{a,3,\{0\}}$ holds and so assuming the ACC  in dim 3 we can prove that $a=0$. Now assume that  $BAB_{\epsilon',3,\{0\}}$ holds for any $\epsilon'>\epsilon$ where $1>\epsilon>0$.

\item Prove $SC_{\epsilon,3}$ in the local case. Moreover prove that the local  $\epsilon$-lc complement indexes can be chosen such that there is a $\tau >0$ s.t. if  $1-\epsilon-\tau \leq b\leq 1-\epsilon$ then $\llcorner (n+1)b\lrcorner \geq n(1-\epsilon)$ for any local $\epsilon$-lc complement index $n$. 

\item Blow up all exceptional divisor $E$ such that $\epsilon\leq a(E,X,0)\leq \epsilon+\tau$ to get $f: Y\longrightarrow X$. Then  $D_{\tau, \{1-\epsilon\}}:=\sum_{i}(1-\epsilon)B_{i}$ where $B_Y=\sum_{i}b_iB_i$ is the crepant pull back boundary. Then run the $D$-LMMP for $D:=-(K_Y+D_{\tau, \{1-\epsilon\}})$. At the end we get $g:Y\dasharrow X_1$ and $X_1\dasharrow S_1$ such that $-(K_{X_1}+D_{\tau, \{1-\epsilon\}})$ is nef and $\equiv 0/S_1$ and $-(K_{S_1}+D_{\tau, \{1-\epsilon\}}).R>0$ for any birational type extremal ray $R$. 

\item There are the following possibilities for the model $S_1$:

\begin{description}

\item[$\diamond$] $\rho(S_1)=1$, $K_{S_1}+D_{\tau, \{1-\epsilon\}}$ is ample and $K_{S_1}+D_{\tau, \{1-\epsilon\}}$ is $\epsilon$-lc. 

\item[$\diamond$] $-(K_{S_1}+D_{\tau, \{1-\epsilon\}}).R=0$ for a fibre type extremal ray $R$ and the log divisor  $K_{S_1}+D_{\tau, \{1-\epsilon\}}$ is $\epsilon$-lc.

\item[$\diamond$]  $-(K_{S_1}+D_{\tau, \{1-\epsilon\}})$ is nef and big and $K_{S_1}+D_{\tau, \{1-\epsilon\}}$ is $\epsilon$-lc.
\end{description}

\item If the first case happens in the division above then delete the boundary, so $(S_1, 0)$ is $\epsilon+\tau$-lc and so the pair is bounded by the assumptions. 

\item \begin{defn}
Let $f: T\longrightarrow Z$ be a contraction and $K_T+B\sim_{\mathbb{R}}0/Z$. Put $D_Z:=\sum_{i}d_iD_i$ where $d_i$ is defined as follows:

 \[1-d_i=\sup \{c|K_T+B+c{f^*D_i} ~\mbox{is lc over the generic point of $D_i$}\}\]

\end{defn}

\item If the second case occurs in the division above then we need the following general conjecture, due to Shokurov [PSh1] and Kawamata [K3], which is useful in many situations:

\begin{conj}[Adjunction]
Let $(T/Z,B)$ be a lc pair of dim $d$ such that $K_T+B\sim_{\mathbb{R}}0/Z$. Define the unique class $M_Z$ up to $\mathbb{R}$-linear equivalence as $K_T+B\sim_{\mathbb{R}}{^*(K_Z+D_Z+M_Z)}$. Then the followings hold:
\begin{description}
\item[Adjunction] We can choose an $M_Z\geq 0$ in its $\mathbb{R}$-linear equivalence class such that $(Z,D_Z+M_Z)$ is lc.
\item[Effective adjunction] Fix $\Gamma_f$. Then there is a constant $I\in \mathbb{N}$ depending only on $d$ and $\Gamma_f$ such that $|IM_Z|$ is a free linear system for an appropriate choice of $M_Z$. In addition the following holds
\[I(K_T+B)\sim {^*I(K_Z+D_Z+M_Z)}.\] 
\end{description} 
\end{conj}

It is expected that the effective adjunction implies the boundedness of $S_1$ under our assumptions. 

\item If the third case occurs in the division above then we need to repeat the pcocess with a bigger $\epsilon$. We have new coefficients in the boundary. Moreover we need to prove that this process stops after a bounded number of steps.

\item If every time the third case happens then at the end we get a pair $(X_r,B_r)$  which is terminal in codim 2 and  $-(K_{X_r}+B_r)$ is nef and big. After some more blow ups and blow downs we may prove that the index of $K_{X_r}+B_r$ is bounded.

\end{enumerate}

\clearpage

%%%%%%%%%%%%%%%%%%%%%%%%%%%%%%%%%%%%%%%%%%%%%%%%%%%%%%%%%%%%%%%%%%%%%%

\section{{List of notation and terminology}}

\begin{tabular*}{12cm}{l l l}

 ${\mathbf{\mathbb{N}}}$  &&  \parbox[t]{10cm} {\emph{The set of natural numbers $\{1,2, \dots\}$.}} \\  \\

 ${\mathbf{\mathbb{R}^+}}$  &&  \parbox[t]{10cm} {\emph{The set of positive real numbers. Similar notation for $\mathbb{Q}$.}} \\ \\

  {\textbf{dim}}  &&  \parbox[t]{10cm} {\emph{Dimension or dimensional.}} \\ \\
  
 {\textbf{WLF}}  && \parbox[t]{10cm}{\emph{Weak log Fano. $(X/Z,B)$ is WLF if $X/Z$ is a projective contraction and $-(K_{X}+B)$ is nef and big/$Z$ and $X$ is $\mathbb{Q}$-factorial.}}  \\ \\

  {\textbf{Pseudo-WLF}} && \parbox[t]{10cm}{\emph{Pseudo weak log Fano/$Z$ i.e. there is a $B$ where  $(X/Z,B)$ is WLF.}}\\ \\
 
  ${\mathbf{\Phi_{sm}}}$  &&  \parbox[t]{10cm} {\emph{The set of standard boundary coefficients i.e. $\{\frac{k-1}{k}\}_{k\in \mathbb{N}}\cup \{1\}$.}} \\ \\

  ${\mathbf{\Gamma_f}}$  &&  \parbox[t]{10cm} {\emph{A finite subset of $[0,1]$.}} \\ \\
 
  ${\mathbf{\mld(\mu,X,B)}}$  && \parbox[t]{10cm}{\emph{The log minimal discrepancy of $(X,B)$ at the centre $\mu$.}} \\ \\

  ${\mathbf{index_{P}(D)}}$  && \parbox[t]{10cm}{\emph{The smallest positive natural number $r$ s.t. $rD$ is a Cartier divisor at $P$.}} \\ \\

  ${\mathbf{WC_{\delta,d,\Gamma}}}$  &&\parbox[t]{10cm} {\emph{The weak conjecture on the boundedness of $\epsilon$-lc complements in dim $d$. See \ref{weak}}} \\ \\

  $\mathbf{SC_{\delta,d}}$  &&\parbox[t]{10cm} {\emph{The strong conjecture on the boundedness of $\epsilon$-lc complements in dim $d$. See \ref{strong}}} \\ \\

  $\mathbf{BAB_{\delta,d,\Gamma}}$  && \parbox[t]{10cm}{\emph{The Alexeev-Borisovs conjecture on the boundedness of $d$-dim $\delta$-lc WLF varieties. See \ref{BAB}}} \\ \\

  ${\mathbf{LT_{d}}}$  && \parbox[t]{10cm}{\emph{The log termination conjecture in dim $d$. See \ref{lt}}} \\ \\

  ${\mathbf{ACC_{d,\Gamma}}}$  &&\parbox[t]{10cm} {\emph{The ACC conjecture on mlds in dim $d$. See \ref{acc}}} \\ \\

\end{tabular*}

%%%%%%%%%%%%%%%%%%%%%%%%%%%%%%%%%%%%%%%%%%%%%%%%%%%%%%%%%%%%%%%%%

\section{{References:}}

\begin{tabular*}{12cm}{l l l}
{\textbf{[A1]}} & & \parbox[t]{11cm} {V. Alexeev; {\emph{Boundedness and $K\sp 2$ for log surfaces.}}  Internat. J. Math.  5  (1994),  no. 6, 779--810.}\\ \\

{\textbf{[A2]}} & & \parbox[t]{11cm} {V. Alexeev; {\emph{Two two dimensional terminations.}} Duke Math. J.  69  (1993),  no. 3, 527--545. }\\ \\

{\textbf{[AM]}}&   &\parbox[t]{11cm} {V. Alexeev, S. Mori; {\emph{Bounding singular surfaces of general type.}} Algebra, arithmetic and geometry with applications (West Lafayette,IN, 2000), 143--174, Springer, Berlin, 2004.}\\ \\

{\textbf{[Am]}} & & \parbox[t]{11cm} {F. Ambro; {\emph{On minimal log discrepancies.}}  Math. Res. Lett.  6  (1999),  no. 5-6, 573--580.}\\ \\

{\textbf{[B]}}&  & \parbox[t]{11cm} {A.A. Borisov;  {\emph{Boundedness of Fano threefolds with log-terminal singularities of given index.}}  J. Math. Sci. Univ. Tokyo  8  (2001), no. 2, 329--342.}\\ \\

{\textbf{[C]}}&  & \parbox[t]{11cm} {A. Corti;  {\emph{Recent results in higher-dimensional birational geometry.}} Current topics in complex algebraic geometry (Berkeley, CA, 1992/93),  35--56, Math. Sci. Res. Inst. Publ., 28, Cambridge Univ. Press, Cambridge, 1995. }\\ \\

{\textbf{[C1]}}&  & \parbox[t]{11cm} {A. Corti;  {\emph{Singularities of linear systems and $3$-fold birational geometry.}} Explicit birational geometry of 3-folds,  259--312, London Math. Soc. Lecture Note Ser., 281, Cambridge Univ. Press, Cambridge, 2000. }\\ \\

{\textbf{[CR]}}&  & \parbox[t]{11cm} {A. Corti, M. Reid;  {\emph{Explicit birational geometry of 3-folds.}} Edited by Alessio Corti and Miles Reid. London Mathematical Society Lecture Note Series, 281. Cambridge University Press, Cambridge, 2000. }\\ \\

{\textbf{[CPR]}}&  & \parbox[t]{11cm} {A. Corti, A. Pukhlikov, M. Reid; {\emph{Fano $3$-fold hypersurfaces.}}  Explicit birational geometry of 3-folds,  175--258, London Math. Soc. Lecture Note Ser., 281, Cambridge Univ. Press, Cambridge, 2000.}\\

\end{tabular*}

\begin{tabular*}{12cm}{l l l}

{\textbf{[H]}}&   & \parbox[t]{11cm} {R. Hartshorne; {\emph{Algebraic geometry.}} Graduate Texts in Mathematics, No. 52. Springer-Verlag, 1977.}\\ \\

{\textbf{[K1]}}&   &\parbox[t]{11cm} {Y. Kawamata; {\emph{Boundedness of $\bold Q$-Fano threefolds.}}  Proceedings of the International Conference on Algebra, Part 3 (Novosibirsk, 1989),  439--445, Contemp. Math., 131, Part 3, Amer. Math. Soc., Providence,  RI, 1992.}\\ \\

{\textbf{[K2]}}&&   \parbox[t]{11cm} {Y. Kawamata; {\emph{Termination of log flips for algebraic $3$-folds.}} Internat. J. Math.  3  (1992),  no. 5, 653--659.}\\ \\

{\textbf{[K3]}}&& \parbox[t]{11cm} {Y. Kawamata; {\emph{  Subadjunction of log canonical divisors for a subvariety of codimension $2$}}.  Birational algebraic geometry (Baltimore, MD, 1996),  79--88, Contemp. Math., 207, Amer. Math. Soc., Providence, RI, 1997.}\\ \\

{\textbf{[K4]}}&& \parbox[t]{11cm} {Y. Kawamata; {\emph{ Subadjunction of log canonical divisors. II}}.  Amer. J. Math.  120  (1998),  no. 5, 893--899. }\\ \\

{\textbf{[K5]}}&& \parbox[t]{11cm} {Y. Kawamata; {\emph{  The number of the minimal models for a $3$-fold of general type is finite}}.  Math. Ann.  276  (1987),  no. 4, 595--598.}\\ \\

{\textbf{[K6]}}&& \parbox[t]{11cm} {Y. Kawamata; {\emph{ Termination of log flips in dimension 4}}. Preprint. It contained a proof which turned to be not correct.}\\ \\

{\textbf{[KMM] }}&& \parbox[t]{11cm} {Y. Kawamata, K. Matsuda, K. Matsuki; {\emph{Introduction to the minimal model problem.}}  Algebraic geometry, Sendai, 1985,  283--360, Adv. Stud. Pure Math.,10, North-Holland, Amsterdam, 1987.}\\ \\

{\textbf{[Ko1] }}&& \parbox[t]{11cm} {J. Kollar; {\emph{Singularities of pairs.}}  Algebraic geometry---Santa Cruz 1995,  221--287, Proc. Sympos. Pure Math., 62, Part 1, Amer. Math. Soc., Providence, RI, 1997.} \\ \\

\end{tabular*}

\begin{tabular*}{12cm}{l l l}

{\textbf{[Ko2] }}&& \parbox[t]{11cm} {J. Kollar; {\emph{ Rational curves on algebraic varieties.}} A Series of Modern Surveys in Mathematics [Results in Mathematics and Related Areas. 3rd Series. A Series of Modern Surveys in Mathematics], 32. Springer-Verlag, Berlin, 1996.} \\ \\

{\textbf{[KD] }}&& \parbox[t]{11cm} {J. Demailly, J. Kollar; {\emph{Semi-continuity of complex singularity exponents and Kähler-Einstein metrics on Fano orbifolds.}} Ann. Sci. École Norm. Sup. (4)  34  (2001),  no. 4, 525--556. } \\ \\

{\textbf{[KM]}}&&   \parbox[t]{11cm} {J. Kollar, S. Mori; {\emph{Birational geometry of algebraic varieties.}} With the collaboration of C. H. Clemens and A. Corti. Translated  from the 1998 Japanese original. Cambridge Tracts in Mathematics, 134. Cambridge University Press, Cambridge, 1998.} \\ \\

{\textbf{[KMMT]}}& & \parbox[t]{11cm} {J. Koll\'ar, Y. Miyaoka, S. Mori, H. Takagi; {\emph{Boundedness of canonical $\bold Q$-Fano 3-folds.}} Proc. Japan Acad. Ser. A Math. Sci.  76  (2000),  no. 5, 73--77.}\\ \\

${\mathbf{[K^+]}}$& & \parbox[t]{11cm} {J. Koll\'ar and others; {\emph{Flips and abundance for algebraic threefolds.}} Papers from the Second Summer Seminar on Algebraic Geometry held at the University of Utah, Salt Lake City, Utah, August 1991. Astérisque No. 211 (1992). Société Mathématique de France, Paris,  1992. pp. 1--258.}\\ \\

{\textbf{[Mc]}}&   & \parbox[t]{11cm} {J. McKernan; {\emph{ Boundedness of log terminal Fano pairs of bounded index.}} ArXiv/math.AG/0205214}\\ \\

{\textbf{[MP]}}&   & \parbox[t]{11cm} {J. McKernan, Yu. Prokhorov;{\emph{ Threefold Thresholds.}} ArXiv/math.AG/0205214}\\ \\

{\textbf{[Pr] }}&& \parbox[t]{11cm} { Yu. Prokhorov; {\emph{ Lectures on complements on log surfaces.}} MSJ  Memoirs, 10. Mathematical Society of Japan, Tokyo, 2001.}\\ \\

{\textbf{[Pr1] }}&& \parbox[t]{11cm} { Yu. Prokhorov; {\emph{Boundedness of exceptional quotient singularities. }}(Russian)  Mat. Zametki  68  (2000),  no. 5, 786--789;  translation in  Math. Notes  68  (2000),  no. 5-6, 664--667}\\ \\

\end{tabular*}

\begin{tabular*}{12cm}{l l l}

{\textbf{[Pr2] }}&& \parbox[t]{11cm} { Yu. Prokhorov; {\emph{Boundedness of nonbirational extremal contractions.}}  Internat. J. Math.  11  (2000),  no. 3, 393--411. }\\ \\

{\textbf{[PrM] }}&& \parbox[t]{11cm} {D. Markushevich, Yu. Prokhorov; {\emph{Exceptional quotient singularities.}}  Amer. J. Math.  121  (1999),  no. 6, 1179--1189. }\\ \\

{\textbf{[PrI] }}&& \parbox[t]{11cm} {V.A. Iskovskikh, Yu. Prokhorov; {\emph{ Fano varieties.  Algebraic geometry, V,  1--247, Encyclopaedia Math. Sci., 47, Springer, Berlin, 1999.}} }\\ \\

{\textbf{[PSh]}}&&   \parbox[t]{11cm} {Yu. Prokhorov; V.V. Shokurov; {\emph{The first fundamental theorem on complements: from global to local.}} (Russian)  Izv. Ross. Akad.  Nauk Ser. Mat.  65  (2001),  no. 6, 99--128;  translation in  Izv. Math.  65  (2001),  no. 6, 1169--1196.}\\ \\

{\textbf{[PSh1]}}&&  \parbox[t]{11cm} {Yu. Prokhorov; V.V. Shokurov; {\emph{Toward the second main theorem on complements: from local to global.}} Preprint 2001.}\\ \\

{\textbf{[R]}}&&    \parbox[t]{11cm} {M. Reid; {\emph{Chapters on algebraic surfaces.}}  Complex algebraic  geometry (Park City, UT, 1993),  3--159, IAS/Park City Math. Ser., 3, Amer. Math. Soc., Providence, RI, 1997.}\\ \\

{\textbf{[R1]}}&&    \parbox[t]{11cm} {M. Reid; {\emph{Update on 3-folds.}}  Proceedings of the International Congress of Mathematicians, Vol. II (Beijing, 2002),  513--524, Higher Ed. Press, Beijing, 2002.  }\\ \\

{\textbf{[R2]}}&&    \parbox[t]{11cm} {M. Reid; {\emph{Twenty-five years of $3$-folds---an old person's view.}} Explicit birational geometry of 3-folds,  313--343, London Math. Soc. Lecture Note Ser., 281, Cambridge Univ. Press, Cambridge, 2000. }\\ \\

{\textbf{[R3]}}&&    \parbox[t]{11cm} {M. Reid; {\emph{Young person's guide to canonical singularities.}} Algebraic geometry, Bowdoin, 1985 (Brunswick, Maine, 1985),  345--414, Proc. Sympos. Pure Math., 46, Part 1, Amer. Math. Soc., Providence, RI, 1987. }\\ \\

{\textbf{[R4]}}&&    \parbox[t]{11cm} {M. Reid; {\emph{The moduli space of $3$-folds with $K=0$ may nevertheless be irreducible.}} Math. Ann.  278  (1987),  no. 1-4, 329--334. }\\ \\

\end{tabular*}

\begin{tabular*}{12cm}{l l l}

{\textbf{[Sh1]}}&&  \parbox[t]{11cm} {V.V. Shokurov; {\emph{Three-dimensional log flips.}} With an appendix in English by Yujiro Kawamata.  Russian  Acad. Sci. Izv. Math.  40  (1993),  no. 1, 95--202.}\\ \\

{\textbf{[Sh2]}}&&  \parbox[t]{11cm} {V.V. Shokurov; {\emph{Complements on surfaces.}} Algebraic geometry,  10.  J. Math. Sci. (New York)  102  (2000),  no. 2, 3876--3932.}\\ \\

{\textbf{[Sh3]}}&&  \parbox[t]{11cm} {V.V. Shokurov; {\emph{Prelimiting flips.}}  Tr. Mat. Inst. Steklova  240  (2003),  Biratsion. Geom. Linein. Sist. Konechno Porozhdennye  Algebry, 82--219;  translation in  Proc. Steklov Inst. Math.  2003, no. 1 (240), 75--213.} \\ \\

{\textbf{[Sh4]}}&&  \parbox[t]{11cm} {V.V. Shokurov; {\emph{Letters of a birationalist  V: Mld's and termination of log flips}}.} \\ \\

{\textbf{[Sh5]}}&&  \parbox[t]{11cm} {V.V. Shokurov; {\emph{$3$-fold log models.}} Algebraic geometry, 4.  J. Math.  Sci.  81  (1996),  no. 3, 2667--2699.}\\ \\

{\textbf{[Sh6]}}&&  \parbox[t]{11cm} {V.V. Shokurov; {\emph{Letters of a bi-rationalist. IV. Geometry of log flips.}}  Algebraic geometry,  313--328, de Gruyter, Berlin, 2002.}\\ \\

{\textbf{[Sh7]}}&&  \parbox[t]{11cm} {V.V. Shokurov; {\emph{A nonvanishing theorem.}} (Russian)  Izv. Akad. Nauk  SSSR Ser. Mat.  49  (1985),  no. 3, 635--651.}\\ \\

{\textbf{[Sh8]}}&&  \parbox[t]{11cm} {V.V. Shokurov; {\emph{ACC in codim 2}}. Preprint.}\\ \\

\end{tabular*}

%%%%%%%%%%%%%%%%%%%%%%%%

\end{document}